\documentclass{article}

\usepackage[english]{babel}
\usepackage{amssymb,amsmath,amsthm}

\theoremstyle{plain}
\newtheorem{proposition}{Proposition}
\newtheorem{theorem}[proposition]{Theorem}
\newtheorem{lemma}[proposition]{Lemma}
\newtheorem{corollary}[proposition]{Corollary}
\theoremstyle{definition}
\newtheorem{definition}[proposition]{Definition}
\newtheorem{example}[proposition]{Example}
\theoremstyle{remark}
\newtheorem*{remark}{Remark}

\usepackage[english]{babel}
\newcommand{\R}{\mathbb{R}}
\newcommand{\im}{\mathrm{Im}\,}         

\newcommand{\Lie}{\mathcal{L}}          
\newcommand{\C}{\mathbb{C}}
\newcommand{\hook}{\lrcorner\,}
\newcommand{\LieG}[1]{\mathrm{#1}}      
\newcommand{\Spin}{\LieG{Spin}}
\newcommand{\SU}{\LieG{SU}}
\newcommand{\SO}{\LieG{SO}}
\newcommand{\GL}{\LieG{GL}}
\newcommand{\SL}{\LieG{SL}}
\newcommand{\Gtwo}{\LieG{G}_2}

\newcommand{\ann}{^\mathrm{o}}

\newcommand{\lie}[1]{\mathfrak{#1}}     
\newcommand{\so}{\lie{so}}
\newcommand{\gl}{\lie{gl}}
\newcommand{\su}{\lie{su}}
\newcommand{\spin}{\lie{spin}}
\newcommand{\dfn}[1]{\emph{#1}}
\newcommand{\id}{\mathrm{Id}}   

\DeclareMathOperator{\tr}{tr}
\DeclareMathOperator{\Sym}{Sym}
\DeclareMathOperator{\ric}{Ric}

\renewcommand{\Re}{\mathfrak{Re}}
\renewcommand{\Im}{\mathfrak{Im}}
\newcommand{\ee}[1]{e^#1\kern-2pt\otimes\kern-1pt e^#1}

\begin{document}
\title{Generalized Killing Spinors in Dimension 5}
\author{Diego Conti\\ Simon Salamon}
\date{29th June 2006}
\maketitle
\begin{abstract}We study the intrinsic geometry of hypersurfaces in Calabi-Yau manifolds of real dimension 6
and, more generally, $\mathrm{SU}(2)$-structures on 5\nobreakdash-\hspace{0pt}manifolds defined by a generalized Killing
spinor. We prove that in the real analytic case, such a 5-manifold can be isometrically embedded as a hypersurface in a
Calabi-Yau manifold in a natural way.  We classify nilmanifolds carrying invariant structures of this type, and present
examples of the associated metrics with holonomy $\mathrm{SU}(3)$.\end{abstract}

\vskip5pt\centerline{\small\textbf{MSC classification}: 53C25; 14J32, 53C29, 53C42, 58A15}\vskip15pt

Let $N$ be a spin manifold, and let $\Sigma_N$ be the complex spinor bundle, which splits as
$\Sigma^+_N\oplus\Sigma^-_N$ in even dimension. It is well known that any oriented hypersurface $\iota\colon M\to N$ is
also spin, and we have $\Sigma_M=\iota^*\Sigma_N$ or $\Sigma_M=\iota^*\Sigma_N^+$ according to whether the dimension of
$N$ is odd or even. Thus, a spinor $\psi_N$ on $N$ (which we assume to lie in $\Sigma_N^+$ if the dimension is even)
induces a spinor $\psi=\iota^*\psi_N$ on $M$. If $\psi_N$ is parallel with respect to the Levi-Civita connection of
$N$, then
\begin{equation}
\label{eqn:GeneralizedKillingSpinor} \nabla_X \psi=\frac{1}{2}A(X)\cdot \psi
 \end{equation}
where $\nabla$ is the covariant derivative with respect to the Levi-Civita connection of $M$, the dot represents
Clifford multiplication and $A$ is a section of the bundle of symmetric endomorphisms of $TM$; in fact, $A$ is the
Weingarten tensor. On a Riemannian spin manifold, spinors $\psi$  satisfying (\ref{eqn:GeneralizedKillingSpinor}) for
some symmetric $A$ are called  generalized Killing spinors \cite{BarGauduchonMoroianu}. Generalized Killing spinors
with $\tr(A)$ constant arise in the study of the Dirac operator, and are called $T$-Killing spinors
\cite{FriedrichKim:TheEinsteinDiracEquation}. For a consistent terminology, we define Killing spinors by the condition
$A=\lambda\,\id$, where $\lambda$ is required to be a \emph{real} constant. If $N$ is the cone on $M$, i.e. the warped
product $M\times_{r} \R^+$, then $\iota^*\psi$ is a Killing spinor.

Any generalized Killing spinor $\psi$ is parallel with respect to a suitable connection; consequently, $\psi$ defines a
$G$-structure consisting of those frames $u$ such that \mbox{$\psi=[u,\psi_0]$} for some fixed $\psi_0$ in $\Sigma_n$,
where $G$ is the stabilizer of $\psi_0$. The intrinsic torsion of this $G$-structure can be identified with $A$. It is
easy to prove that the $G$-structures defined by a generalized Killing spinor are cocalibrated $\Gtwo$-structures in
dimension 7 and half-flat $\SU(3)$-structures in dimension 6. The statement that a parallel spinor restricts to a
generalized Killing spinor is therefore a generalization of the following: a hypersurface in an 8\nobreakdash-manifold
with holonomy contained in $\Spin(7)$ (resp.\ a 7-manifold with holonomy contained in $\Gtwo$) inherits a natural
cocalibrated $\Gtwo$-structure (resp.\ half-flat $\SU(3)$\nobreakdash-structure). In the present article, we study the
$\SU(2)$-structures in dimension 5 defined by a generalized Killing spinor; we call these structures \emph{hypo}. By
the above discussion, it is clear that hypersurfaces inside 6-manifolds with holonomy $\SU(3)$ inherit a natural hypo
structure. The word `hypo' reflects the fact that such structures are however \emph{under}-defined in senses that will
become clearer during the course of the paper.

As Killing spinors in dimension 5 correspond to Einstein-Sasaki structures \cite{FriedrichKath}, hypo geometry is a
generalization of Einstein-Sasaki geometry. However, hypo geometry is not closely related to other generalizations like
contact metric structures, or Sasaki structures. Indeed, hypo structures are not necessarily contact, and in Section
\ref{sec:Sasaki} we show that, locally, Sasaki manifolds have a compatible hypo structure if and only if they are $\alpha$-Einstein. On
the other hand, just like half-flat is much weaker than nearly-K\"ahler, so is hypo much weaker than Einstein-Sasaki, or even $\alpha$-Einstein-Sasaki. In fact, the first Betti number $b_1$ of a compact  Sasaki manifold is even, a condition which need not be satisfied by compact hypo manifolds (for instance, the product of a hyperk\"ahler 4-manifold with $S^1$ is hypo; for
less trivial examples, see Section \ref{sec:nilmanifolds}).

It is natural to ask whether any spin Riemannian manifold $M$ with a generalized Killing spinor can be embedded as a
hypersurface in some $N$ so that the above construction gives back the starting metric and spinor on $M$. If so, we say
that $M$ has the \textit{embedding property}. Notice that we are not only requiring the embedding to be isometric, but
also that the Weingarten tensor coincide with the intrinsic torsion $A$. The embedding property is known to hold in
dimension 2 \cite{Friedrich:OnTheSpinorRepresentationOfSurfaces} and in the case in which $A$ satisfies the Codazzi
equation \cite{BarGauduchonMoroianu} (see also \cite{Morel} for the case in which $A$ is parallel). In dimension 6  and 7, Hitchin has translated the embedding property into a
problem of existence of integral lines for a vector field on an infinite-dimensional space \cite{Hitchin:StableForms}.
This approach can be adapted to the hypo situation, although the number of differential forms required in the
definition (one 1-form and three 2-forms) complicates matters \cite{thesis} and raises a number of questions that will
be addressed elsewhere.

In Section \ref{sec:EmbeddingProperty} we use Cartan-K\"ahler theory to prove the following  theorem: every
real analytic 5-manifold with a real analytic hypo structure has the embedding property. To establish this
result, we use the well-known fact that the corresponding differential system is involutive
\cite{Bryant:Calibrated}, and with methods closely based on \cite{Bryant:Calibrated} we construct a regular
flag at each point whose 5-dimensional element is tangent to the 5-dimensional integral manifold given by the
hypo structure; a standard argument of Cartan-K\"ahler theory completes the proof. An analogous result is
expected to hold for hypersurfaces of spaces with Ricci-flat holonomy group, including $\SU(n)$. Note that in
dimension 2, the Codazzi equation (imposing that $\nabla A$ be totally symmetric) is automatically satisfied,
but this is not true in dimension 5.

In Section~\ref{sec:nilmanifolds}, we give a complete list of the 5-dimensional nilmanifolds which admit an invariant
hypo structure.  The analogous classification problem in dimension 6 is still open, although of the 34 isomorphism
classes of 6-dimensional nilmanifolds, 12 are known to admit invariant half-flat structures
\cite{ChiossiFino,ChiossiSwann}, and we have been able to produce 11 more. All of the resulting compact 5-dimensional
examples satisfy the hypotheses of the embedding theorem, and allow one to construct a host of explicit Ricci-flat
(albeit incomplete) metrics with holonomy group equal to $\SU(3)$.

 We illustrate the importance of these examples in Section~\ref{sec:examples} by showing first that they do not in
general satisfy the Codazzi equation. Finally, we explain that for an appropriate choice of nilpotent Lie algebra, the
construction gives rise to metrics with holonomy $\SU(3)$ with  the following \textit{tri-Lagrangian}
property: there is an $S^1$ family of Lagrangian submanifolds passing through each point, including three that are
special Lagrangian for any fixed phase.

\section{$\SU(2)$-structures and  hypersurfaces} In this section we show how $\SU(2)$-structures, and hypo
structures in particular, can be defined using differential forms, and in this language we carry out a construction
described in the introduction.  Namely, we show that $\SU(2)$\nobreakdash-structures arise naturally on hypersurfaces of
6-manifolds endowed with an $\SU(3)$-structure, and that if the $\SU(3)$-structure is integrable, the corresponding
$\SU(2)$-structure satisfies the hypo condition which is characterized by Definition~\ref{dfn:hypo}. A particular case
is when the 6-manifold is a cone on the 5-manifold; in that case, the metric induced on the latter is Einstein-Sasaki
\cite{BoyerGalicki}.  In fact, hypo structures generalize Einstein-Sasaki structures: the former are defined by a
generalized Killing spinor, and the latter by a Killing spinor, as will be illustrated in Section~\ref{section:IntrinsicTorsion}.  

Let $M$ be a 5-manifold.  An $\SU(2)$-structure on $M$ is an $\SU(2)$-reduction $P$
of the frame bundle $F$ on $M$, or equivalently a section $s$ of the bundle $F/\SU(2)$; more in the spirit of special
geometries, we have the following characterization:

\begin{proposition}
\label{prop:SU2forms} $\SU(2)$-structures on a 5-manifold are in one-to-one correspondence with quadruplets
$(\alpha,\omega_1,\omega_2,\omega_3)$, where $\alpha$ is a 1-form and $\omega_i$ are 2\nobreakdash-forms, satisfying:
\begin{equation}
\label{eqn:omegaiorthogonal} \omega_i\wedge\omega_j=\delta_{ij}\upsilon
\end{equation}
for some  4-form $\upsilon$ with $\upsilon\wedge\alpha\neq 0$, and
 \begin{equation}
 \label{eqn:omegaioriented}X\hook\omega_1=Y\hook\omega_2 \implies
\omega_3(X,Y)\geq 0\;.
\end{equation}
 Equivalently, an $\SU(2)$-structure can be defined by a 1-form $\alpha$, a 2-form $\omega_1$ and a
complex 2-form $\Phi$, corresponding to $\omega_2+i\omega_3$, such that
\begin{align*}
\alpha\wedge\omega_1^2&\neq 0 &\omega_1\wedge\Phi&=0\\
\Phi^2&=0 & 2\omega_1^2&=\Phi\wedge\overline\Phi
\end{align*}
and $\Phi$ is $(2,0)$ with respect to $\omega_1$.
\end{proposition}
\begin{remark} If we start with an $\SO(5)$-structure, we can understand a reduction to $\SU(2)$ as follows. The form
$\alpha$ defines a splitting $\R^5=\R\oplus\R^4$; a metric and an orientation is induced on $\R^4$. The eigenspace
decomposition relative to the Hodge star gives $\Lambda^2(\R^4)^*=\Lambda^2_+\oplus\Lambda^2_-$, which corresponds to
writing $\SO(4)$ as $\SU(2)_+\SU(2)_-$. The choice of a basis $\omega_1,\omega_2,\omega_3$ of $\Lambda^2_+$,
corresponding to (\ref{eqn:omegaiorthogonal}), reduces then $\SO(4)$ to $\SU(2)_-$. Since $\Lambda^2_+$ has a natural
orientation, one can always assume that $\omega_1,\omega_2,\omega_3$ be a positively oriented, orthogonal basis of this
space; this assumption
corresponds to (\ref{eqn:omegaioriented}).
\end{remark}
Bearing this construction in mind, we shall sometimes refer to the structure group as $\SU(2)_-$ rather than $\SU(2)$.
By $(\alpha,\omega_i)$ we shall mean a quadruplet $(\alpha,\omega_1,\omega_2,\omega_3)$ satisfying Proposition
\ref{prop:SU2forms}, and for a form $\omega$ on a manifold $M$, we define
\[\omega\ann=\{X\in TM\mid X\hook\omega=0\}\;.\]
Before proving Proposition \ref{prop:SU2forms}, it is convenient to prove the following:
\begin{proposition}
\label{prop:dualofalpha} For $(\alpha,\omega_i)$ as above, we have $\omega_i\ann=\upsilon\ann$, $i=1,2,3$.
\end{proposition}
\begin{proof}
Define $\beta_i$ and $\gamma_i$ by the condition that $\omega_i=\alpha\wedge\beta_i+\gamma_i$, with $\beta_i\ann$ and
$\gamma_i\ann$ containing $\upsilon\ann$. Then by (\ref{eqn:omegaiorthogonal}),
\[2\,\alpha\wedge\beta_i\wedge\gamma_i+\gamma_i^2=\upsilon\;.\]
Take a non-zero $X$ in $\upsilon\ann$; then
\[0=X\hook(\alpha\wedge\beta_i\wedge\gamma_i)=(X\hook\alpha)\beta_i\wedge\gamma_i\;,\] and so $\beta_i\wedge\gamma_i=0$. On the other hand,
$0\neq\alpha\wedge\upsilon=\alpha\wedge\gamma_i^2$ shows that $\gamma_i$ is non-degenerate on $\alpha\ann$, and
therefore $\beta_i=0$. Thus $\omega_i\ann\supseteq\upsilon\ann$; the opposite inclusion follows from
$\upsilon=\omega_i^2$.
\end{proof}

\begin{corollary}
\label{cor:referenceformsSU2} Given $(\alpha,\omega_i)$ as above, one can always find a local basis $e^1,\dotsc,e^5$ of forms
on $M$ such that
\begin{equation}
\label{eqn:referenceformsSU2}\left\{
\begin{aligned}\alpha&=e^5&\qquad\omega_1&=e^{12}+e^{34}\\
\omega_2&=e^{13}+e^{42}&\qquad\omega_3&=e^{14}+e^{23}
\end{aligned}\right.
\end{equation}
Moreover, one may require that (for example) $e^1$ equal a fixed unit form orthogonal to $\alpha$.
\end{corollary}
\begin{proof} Write $TM=\alpha\ann\oplus\upsilon\ann$; the restriction of $\omega_i$, $\upsilon$ to $\alpha\ann$ satisfy
(\ref{eqn:omegaiorthogonal}) and (\ref{eqn:omegaioriented}); the statement then follows from its 4-dimensional
analogue.
\end{proof}
Here and in the sequel, $e^{12}$ is short for $e^1\wedge e^2$, and so on. A consequence of the corollary itself is that
a global nowhere vanishing $1$\nobreakdash-form in $\alpha^\perp$ exists only if $M$ is parallelizable; in general,
(\ref{eqn:referenceformsSU2}) can only be used locally.

\begin{proof}[Proof of Proposition \ref{prop:SU2forms}] It is sufficient to show that if
$e^1,\dotsc,e^5$ is the standard basis of $(\R^5)^*$ and $(\alpha,\omega_i)$ are as in (\ref{eqn:referenceformsSU2}),
the stabilizers of $\alpha$, $\omega_1$, $\omega_2$ and $\omega_3$ have intersection $\SU(2)$. In fact, if $A\in
\GL(5,\R)$ preserves these forms, it must preserve the splitting $\R^5=\alpha\ann\oplus\upsilon\ann$, so that
\[A=\begin{pmatrix}B&0\\0&1\end{pmatrix},\;B\in \GL(4,\R)\;.\] On the other hand, such an $A$ preserves
$\omega_1,\omega_2,\omega_3$ if and only if $B$ preserves the standard hyperk\"ahler structure on $\R^4$,
i.e. $B$ lies in $\LieG{Sp}(1)=\SU(2)$.
\end{proof}

\begin{remark}
One has to fix a choice of reference forms on $\R^5$ in order to actually identify an $\SU(2)$-structure with a
quadruplet of forms $(\alpha,\omega_i)$, or a triplet $(\alpha,\omega_1,\Phi)$; we shall henceforth use
(\ref{eqn:referenceformsSU2}) to do so, and associate to a frame $u\colon\R^5\to T_xM$ forms $(\alpha,\omega_i)$ such
that $u^*\alpha=e^5$, and so on.
\end{remark}

\begin{proposition}
\label{prop:SU2fromSU3} Let $\iota\colon M\to N$ be an immersion of an oriented 5\nobreakdash-manifold into a
6-manifold. Then an $\SU(3)$-structure on $N$ defines an $\SU(2)$\nobreakdash-structure on $M$ in a natural way.
Conversely, an $\SU(2)$-structure on $M$ defines an $\SU(3)$\nobreakdash-structure on $M\times\R$ in a natural way.
\end{proposition}
\begin{proof}
The $\SU(3)$-structure on $N$ defines a non-degenerate $2$-form $\omega$ and a complex $3$-form $\Psi$ with stabilizer
$\SL(3,\C)$. Since both $M$ and $N$ are oriented, the normal bundle to $M$ has a canonical unit section, which using
the metric lifts to a section $V$ of $\iota^*TN$. Define forms on $M$ by
\[\alpha=-V\hook\omega\;,\quad \Phi=-iV\hook\Psi\;,\quad \omega_1=\iota^*\omega\;.\] Choose a local basis of $1$-forms on $N$
such that $V$ is dual to $e^6$ and
\begin{align}
\label{eqn:referenceformsSU3} \omega&=e^{12}+e^{34}+e^{56}\;, & \Psi=(e^1+ie^2)\wedge(e^3+ie^4)\wedge(e^5+ie^6)\;;
\end{align}
then $\iota^*e^1,\dotsc,\iota^*e^5$ satisfy (\ref{eqn:referenceformsSU2}).

Vice versa, given an $\SU(2)$-structure on $M$, an $\SU(3)$-structure on $M\times\R$ is defined by $(\omega,\Psi)$
given by
\begin{align}
\label{eqn:SU3fromSU2} \omega&=\omega_1+\alpha\wedge dt\;, & \Psi=\Phi\wedge(\alpha+idt)\;,
\end{align}
where $t$ is a coordinate on $\R$.\end{proof}

We are now ready to introduce hypo structures:
\begin{definition}
\label{dfn:hypo} The $\SU(2)$-structure determined by $(\alpha,\omega_i)$ is called hypo if \begin{align}
\label{eqn:hypo} d\omega_1&=0\;,& d(\alpha\wedge\omega_2)&=0\;,&d(\alpha\wedge\omega_3)&=0\;.
\end{align}
\end{definition}
\begin{remark}
If $(\alpha,\omega_i)$ satisfies \eqref{eqn:hypo}, then the $\SU(2)$-structures obtained rotating $\omega_2$ and
$\omega_3$ also satisfy \eqref{eqn:hypo}; moreover, they induce the same metric. This is akin to the case of integrable
(i.e. Calabi-Yau) $\SU(n)$-structures on $2n$-dimensional manifolds, where multiplying the holomorphic $n$-form by
$e^{i\theta}$, with $\theta$ a constant, gives a different integrable structure corresponding to the same metric.
\end{remark}
\begin{proposition}
\label{prop:HypoHypersurfaces} Let $\iota\colon M\to N$ be an immersion of an oriented 5-manifold into a 6-manifold
with an integrable $\SU(3)$-structure. Then the $\SU(2)$-structure induced on $M$ is hypo.
\end{proposition}
\begin{proof}
From \eqref{eqn:referenceformsSU2} and \eqref{eqn:referenceformsSU3} it follows that
$\iota^*\Psi=\Phi\wedge\alpha$; recall also that $\omega_1=\iota^*\omega$. Since $\Psi$ and $\omega$ are
closed, and $\iota^*$ commutes with $d$, $M$ is hypo.
\end{proof}
By construction, Equations \ref{eqn:hypo} are exactly the conditions one obtains on the $\SU(2)$\nobreakdash-structure
induced on a hypersurface in a 6-dimensional manifold with a parallel spinor; it is therefore not surprising that these
structures are characterized by the existence of a generalized Killing spinor, as we shall prove in the next section.

\section{Spinors and intrinsic torsion}
In this section we use generalized Killing spinors to study the intrinsic torsion of hypo structures. \\
Fix a
5-manifold $M$. \label{section:IntrinsicTorsion} An $\SU(2)$-structure on $M$ induces a spin structure on $M$; this is
because the sequence of inclusions
\begin{align*}
\SU(2)_-&<\SO(4)<\SO(5) \intertext{lifts to a sequence} \SU(2)_-&<\Spin(4)<\Spin(5)\;,
\end{align*}
 as
$\Spin(4)=\SU(2)_+\times\SU(2)_-$. One can therefore extend the $\SU(2)$-structure $P$ to
\[P_{\Spin(5)}=P\times_{\SU(2)}\Spin(5)\;.\] More, the spinor bundle is $P\times_{\SU(2)}\Sigma$, where $\Sigma\cong\C^4$ is the
complex spinor representation, and $\Spin(5)$ acts transitively on the sphere in $\Sigma$, with stabilizer (conjugate
to) $\SU(2)$. We shall fix a unit $\psi_0\in\Sigma$ whose stabilizer is $\SU(2)_-$. The conclusion is that
$\SU(2)$-structures $P$ on $M$ are in one-to-one correspondence with pairs $(P_{\Spin(5)},\psi)$, where $P_{\Spin(5)}$
is a spin structure on $M$ and $\psi$ is a unit spinor; explicitly, one has
\[\psi=[s,\psi_0]\] for every local section $s$ of $P$.

Now fix an $\SU(2)$-structure $P$ on $M$; let $\psi$ be the defining spinor and $(\alpha,\omega_i)$ the defining forms.
In this section we shall consider both abstract $\SU(2)$\nobreakdash-structures and $\SU(2)$-structures induced by an
immersion
$\iota\colon M\to N$, where $N$ has an $\SU(3)$-structure; in the latter case, we shall assume that this structure is integrable.
Let $\su (2)^\perp$ be the orthogonal complement of $\su(2)_-$ in $\so(5)$ and set
\[\mathcal{T}=T^*\otimes\su(2)^\perp\;,\]
where $T=\R^5$ as an $\SU(2)$-module. Then $\mathcal{T}$ is naturally a (left) $\SU(2)$-module, and the intrinsic
torsion of $P$ is an $\SU(2)$-equivariant map $\Theta\colon P\to\mathcal{T}$. The intrinsic torsion is determined by
$(\nabla\alpha,\nabla\omega_i)$, where $\nabla$ is the covariant derivative relative to the Levi-Civita connection (see
\cite{Salamon:Redbook}, p. 22); we shall now state the analogous result
for spinors, which is proved in the same way.

It is well known that $\Sigma$ has an inner product preserved by $\Spin(5)$. The infinitesimal action of $\Spin(5)$ on
$\psi_0$ gives a map \[\rho_*\colon\spin(5)\to T_{\psi_0}\Sigma=\Sigma\] with kernel $\su(2)$. Using the fact that
$\rho$ preserves some inner product, we see that upon restricting to $\su(2)^\perp$ we get an inclusion
\[\rho_*\colon\su(2)^\perp\to \psi_0^\perp\;.\]
\begin{proposition}
\label{prop:TorsionInTermsOfSpinors} The intrinsic torsion of $P$ is determined by
\[(\id\otimes\rho_*)\circ\Theta=-\nabla\psi\;,\]
where $\nabla\psi$ is viewed as an equivariant map from $P$ to $T^*\otimes\Sigma$.
\end{proposition}
Now consider the injective map \[C\colon T^*\otimes T\to T^*\otimes\Sigma\] mapping $\eta\otimes X$ to $\eta\otimes
(X\cdot\psi_0)$. Define a subspace of $\mathcal{T}$ by
\[\mathcal{T}^K =(\id\otimes \rho_*)^{-1}(C(\Sym(T)))\;,\] where $\Sym(T)$ is the space of symmetric endomorphisms
of $T$. Since $\id\otimes\rho_*$ is injective on $\mathcal{T}$,
\begin{equation}
\label{eqn:GeneralizedKillingTorsion} \mathcal{T}^K\cong\Sym(T)\;.
\end{equation}
 From \eqref{eqn:GeneralizedKillingSpinor}, we immediately obtain:
\begin{corollary}
\label{cor:GeneralizedKillingTorsion} The spinor $\psi$ is generalized Killing if and only if the intrinsic torsion of
$P$ takes values in $\mathcal{T}^K$. If $P$ is induced by an immersion of $M$ in $N$, under
\eqref{eqn:GeneralizedKillingTorsion} the intrinsic torsion is identified with minus one half the Weingarten tensor.
\end{corollary}
\begin{remark}
All of the above works in any dimension $n$, substituting a suitable Lie group $G$ for $\SU(2)$. However, in dimension
5 $\rho_*$ is an isomorphism because \[\dim\su(2)^\perp=7=\dim\psi_0^\perp\;;\] this only holds if $n$ is not too
large, because the real dimension of $\Sigma$ is $2^{[n/2]+1}$ whereas the dimension of $\Spin(n)$ is $n(n+1)/2$.
\end{remark}

The following is proved independently of the above discussion:
\begin{proposition}
\label{prop:torsion} As an $\SU(2)_-$-module, $\mathcal{T}$ decomposes into irreducible components as follows:
\[\mathcal{T}=7\R\oplus 4\Lambda^1(\R^4)^*\oplus 4\Lambda^2_-(\R^4)^*\;.\]
Write
\begin{align*}
d\alpha&=\alpha\wedge\beta+\sum_{i=1}^3 f^i \omega_i+\omega^-\;,\\
d\omega_i&=\gamma_i\wedge\omega_i+\sum_{j=1}^3 f_i^j\alpha\wedge\omega_j+\alpha\wedge\sigma_i^-.
\end{align*}
Then $f_i^j=\lambda\delta_{ij}+g_i^j$ where $g_i^j=-g_j^i$, and according to the above splitting
\[\Theta(u)=\left((f^i,\lambda,g_i^j),(u^*\beta,u^*\gamma_i),(u^*\omega^-,u^*\sigma_i^-)\right).\]
\end{proposition}
We can rewrite \eqref{eqn:GeneralizedKillingTorsion} as
 \begin{equation}\label{eqn:GeneralizedKillingTorsion2}\mathcal{T}^K\cong 2\R\oplus\Lambda^1(\R^4)^*\oplus 3\Lambda^2_-(\R^4)^*\;;
 \end{equation}
we shall now prove that hypo structures are the $\SU(2)$-structures whose intrinsic torsion takes values in this space.
\begin{proposition}
\label{prop:HypoTorsion} The following are equivalent:
\begin{enumerate}
\item $P$ is hypo.
\item The intrinsic torsion of $P$ has the form
\[\Theta(u)=\left((f^1,0,0,g_2^3),(u^*\beta,0,u^*\gamma_2,u^*\gamma_3),(u^*\omega^-,0,u^*\sigma_2^-,u^*\sigma_3^-)\right)\]
and $ \gamma_2=\beta=\gamma_3$.
\item The spinor $\psi$ is generalized Killing.
\end{enumerate}
\end{proposition}
\begin{proof}
Condition 1 and Condition 2 are equivalent by Proposition \ref{prop:torsion}. Now assume Condition 3 holds; then
setting $A^k(X)=e^k(A(X))$, the Levi-Civita connection form restricted to $P$ can be written as
\begin{equation}
 \label{eqn:HypoConnectionForm}
\begin{pmatrix}
0               &   \omega_{12}     &   \omega_{13} &   \omega_{14} &   -A^2\\
-\omega_{12}    &   0               &   -\omega_{14}&   \omega_{13} &   A^1\\
-\omega_{13}    &   \omega_{14}     &   0           &   -\omega_{12}+A^5&  -A^4\\
-\omega_{14}    &   -\omega_{13}    &\omega_{12}-A^5&  0           &   A^3\\
A^2            &   -A^1           &   A^4        &   -A^3       &   0
\end{pmatrix}\;,
\end{equation}
which reduces to the formula in \cite{FriedrichKath} when $A$ is a scalar multiple of the identity and is proved in the same way. A
straightforward calculation shows that Definition \ref{dfn:hypo} is satisfied; thus Condition 3 implies Condition 1.

On the other hand by Corollary \ref{cor:GeneralizedKillingTorsion}, $\psi$ is a generalized Killing spinor if and only
if the intrinsic torsion lies in $\mathcal{T}^K$, which by \eqref{eqn:GeneralizedKillingTorsion2} is isomorphic to the module where the
intrinsic torsion of a hypo structure takes values. Since Condition 3 implies Condition 2, these isomorphic submodules
of $\mathcal{T}$ coincide and Condition 2 implies Condition 3.
\end{proof}

If $M$ is simply connected, Einstein-Sasaki metrics on $M$ are characterized by the existence of a Killing
spinor with $A=\pm\id$  \cite{FriedrichKath}. Therefore, simply connected Einstein-Sasaki 5-manifolds admit a
hypo structure compatible with the metric. On the other hand, Einstein-Sasaki metrics on $M$ are also
characterized by the condition that the conical metric on $M\times\R_+$ has holonomy contained in $\SU(3)$
(see e.g. \cite{BoyerGalicki}). In order to translate this condition in our language, consider the conical
$\SU(3)$\nobreakdash-structure on $N=M\times\R_+$ induced by $P$, defined by
\begin{align}
\label{eqn:conicalSU3fromSU2} \omega&=t^2\omega_1+t\alpha\wedge dt\;, & \Psi=t^2\Phi\wedge(t\alpha+idt)\;.
\end{align}

\begin{proposition}
\label{prop:ES} The following are equivalent:
\begin{enumerate} \item The conical $\SU(3)$-structure on $M\times\R_+$ induced by $P$ is integrable.  \item The
intrinsic torsion of $P$ has the form \[\Theta(u)=\left((-2,0,0,3),(0,0,0,0),(0,0,0,0)\right)\;.\] \item The spinor
$\psi$ is Killing with $A=-\id$.
\end{enumerate}
\end{proposition}
\begin{proof} From \eqref{eqn:conicalSU3fromSU2},
it immediately follows that Condition 1 is equivalent to
\begin{align} \label{eqn:EinsteinSasaki}
d\alpha&=-2\,\omega_1\;, & d\Phi=-3i\,\alpha\wedge\Phi\;.
\end{align}
On the other hand, by Proposition
\ref{prop:torsion}, \eqref{eqn:EinsteinSasaki} is equivalent to Condition 2. Finally, one can use
\eqref{eqn:HypoConnectionForm} to prove that Condition 3 is also equivalent to \eqref{eqn:EinsteinSasaki}.
\end{proof}

\section{Sasaki structures}\label{sec:Sasaki}
So far, we have seen that hypo geometry is a generalization of Einstein-Sasaki geometry; more generally, it follows from \cite{FriedrichKim:TheEinsteinDiracEquation} that $\alpha$-Einstein-Sasaki manifolds are hypo. In this section we argue that in some sense this is the greatest extent to which hypo geometry can be related to Sasaki geometry.

An \dfn{almost contact metric structure} on a manifold $M$ of dimension $5$ is a $\LieG{U}(2)$-structure on $M$. As
\[\LieG{U}(2)=\LieG{Sp}(2,\R)\cap\SO(5)\;,\] we shall think of an almost contact metric structure as a triple
$(g,\alpha,\omega_1)$, where $g$ is a Riemannian metric, $\alpha$ is a unit 1-form and $\omega_1$ is a unit 2-form.

To an almost contact metric structure, much like to an almost-hermitian structure, one can associate the \dfn{Nijenhuis
tensor}, which is a tensor of type $(2,1)$. Since $\LieG{U}(2)$ is a subgroup of $\LieG{U}(3)$, an almost contact
metric structure on $M$ defines an almost-hermitian structure on the product $M\times\R$; then the Nijenhuis tensors of
$M$ and $M\times\R$ can be identified. If we define a $TM$-valued 1-form $J$ by
\[g(J(X),\cdot)=X\hook\omega\;,\] denoting by $\xi$ the vector field dual to $\alpha$, $N$ is characterized by
\begin{multline}
\label{eqn:Nijenhuis}
N(X,Y)=(\nabla_{JX}J)Y-(\nabla_{JY}J)X+(\nabla_{X}J)(JY)-(\nabla_{Y}J)(JX)+\\
-\alpha(Y)\nabla_X\xi+\alpha(X)\nabla_Y\xi \end{multline}
An almost contact
metric structure is \dfn{quasi-Sasakian} if $N$ vanishes and $\omega$ is closed; it is
\emph{contact} if $d\alpha=-2\,\omega_1$, i.e. \[d\alpha(X,Y)=g(X,J Y)\;.\] A quasi-Sasakian structure is \dfn{Sasaki}
if it is contact.

\begin{proposition}
\label{prop:Nijenhuis} The Nijenhuis tensor of a hypo manifold is
\[N(X,Y)=\alpha(X)\left(J\bigl(A(Y)\bigr)-A(JY)\right)-\alpha(Y)\left(J\bigl(A(X)\bigr)-A(JX)\right)\;.\]
\end{proposition}
\begin{proof}
Omitting summation over $i=1,\dotsc,4$, we have
\begin{align*}
g((\nabla J)Y,Z)&=-2\, e^{i5}(Y,Z)A^i\;,& \nabla\xi&=A^i\otimes Je_i\;;
\end{align*}
writing $A(X,Y)$ for $g(A(X),Y)$, Equation \ref{eqn:Nijenhuis} yields
\begin{multline*}
g(N(X,Y),Z)=2\bigl[-A^i(JX)e^{i5}(Y,Z)+A^i(JY)e^{i5}(X,Z)-A^i(X)e^{i5}(JY,Z)+\\
+A^i(Y)e^{i5}(JX,Z)\bigr]+\alpha(Y)A(X,JZ)-\alpha(X)A(Y,JZ)
\end{multline*}
For $X,Y,Z\in\alpha\ann$:
\begin{align*}
g(N(\xi,Y),Z)&=-A(JY,Z)-A(Y,JZ)=-A(JY,Z)+g(J(A(Y)),Z)\\
g(N(\xi,Y),\xi)&=-A(\xi,JY)\\
g(N(X,Y),\xi)&=-A(JX,Y)+A(JY,X)-A(X,JY)+A(Y,JX)
\end{align*}
proving the statement by the symmetry of $A$.
\end{proof}
\begin{lemma}
\label{lemma:hypoKcontact} A hypo structure  is contact if and only if
\[\omega_1(X,AY)+\omega_1(AX,Y)=-2\,\omega_1(X,Y)\quad\forall\, X,Y\;.\]
\end{lemma}
\begin{proof}
By \eqref{eqn:HypoConnectionForm}, on a hypo manifold
 \begin{equation*}d\alpha(X,Y)=\left(\sum_{i=1}^4A^i\wedge
Je^i\right)(X,Y)=\omega_1(X,AY)+\omega_1(AX,Y)\;.\qedhere\end{equation*}
\end{proof}
\begin{remark}
For $M$ a hypersurface in an $\SU(3)$-holonomy  6-manifold, we recover the condition on the Weingarten tensor
characterizing contact hypersurfaces of a K\"ahler manifold found by Okumura \cite{Blair}.
\end{remark}
In five dimensions, a Sasaki structure is called $\alpha$-Einstein if 
\[\ric=(4-2a)\id +2\alpha\otimes\xi\]
for some function $a$, which must then be  constant. We can now characterize hypo $\SU(2)$-structures which are reductions of Sasaki or quasi-Sasakian
$\LieG{U}(2)$\nobreakdash-structures:
\begin{theorem}
\label{thm:HypoQuasiSasakian} A hypo structure is quasi-Sasakian if and only if $A$ commutes with $J$. A hypo structure
is Sasaki if and only if $A=-\id+a \,\alpha\otimes\xi$ for some constant $a$; a Sasaki $\SU(2)$-structure is hypo if and only if it is $\alpha$-Einstein.
\end{theorem}
\begin{proof}
To prove the first statement, suppose $N(\xi,Y)=0$ for all $Y$ in $\alpha\ann$; by Proposition \ref{prop:Nijenhuis}, it
follows that \[J(A(Y))=A(JY)\;.\] From this very equation it follows that $A(JY,\xi)=0$ for all $JY$, i.e.
\[J(A(\xi))=0=A(J\xi)\;;\] so $A$ commutes with $J$. Conversely, if $A$ commutes with $J$  Proposition \ref{prop:Nijenhuis}
implies that $N=0$.

 By the first part of the theorem and Lemma \ref{lemma:hypoKcontact}, a hypo structure is Sasaki if and only if
\begin{equation}
\label{eqn:HypoSasaki}\left\{\begin{aligned} A(JX)&=J(A(X)) & \forall X&\in TM\\ -2g(JX,Y)&=g(JX,AY)+g(JAX,Y) & \forall
X,Y&\in TM
\end{aligned}\right.\end{equation}
Since $A$ is symmetric, the general solution of \eqref{eqn:HypoSasaki} is
\[A=-\id+a \,\alpha\otimes\xi\] where $a$ is a function. We must prove that if a hypo structure is  Sasaki then $a$ is constant.

For $X$ in $\alpha\ann$, we have
\begin{equation}
\label{eqn:cor:HypoQuasiSasakian1}
\begin{split}
R(X,\xi)\psi&=\nabla_X\nabla_\xi\psi-\nabla_\xi\nabla_X\psi-\nabla_{[X,\xi]}\psi=\\
&=\frac{1}{2}\left((Xa)\xi+a\nabla_X\xi+(1-a)\,\xi\cdot X\right)\cdot\psi\;,\end{split}
\end{equation}
and on the other hand, for generic $X$ (see e.g. \cite{Seminarbericht}):
\[\frac{1}{2}\ric(X)\cdot\psi=\sum_{i=1}^5 e_i\cdot R(e_i,X)\psi\;,\]
where $e_i$ is a local  orthonormal basis, $e_5=\xi$. Using \eqref{eqn:cor:HypoQuasiSasakian1}, we find
\[\ric(\xi)\cdot\psi=\sum_{i=1}^4\bigl((\Lie_{e_i} a)e_i\cdot\xi+a\,e_i\cdot\nabla_{e_i}\xi+(1-a)\,\xi\bigr)\cdot\psi\;.\]
Since on a Sasaki manifold $\nabla_X\alpha=-X\hook\omega_1$, the middle summand acts on $\psi$ as multiplication by $4ai$. On the other hand, every Sasaki 5-manifold satisfies
$\ric(\xi)=4\xi$ \cite{Blair}, so the first summand has to vanish. A similar calculation shows that for $X$ in $\alpha\ann$,
\[\ric(X)\cdot\psi=(Xa+(4-2a)X)\cdot\psi\;;\]
since we already know that $Xa=0$, this implies that the hypo structure is $\alpha$\nobreakdash-Einstein, and in particular $a$ is constant.

 The last statement is now a consequence of the characterization of $\alpha$-Einstein-Sasaki structures in terms of quasi-Sasakian Killing 
spinors \cite{FriedrichKim:TheEinsteinDiracEquation}.
\end{proof}
\section{The embedding property}\label{sec:EmbeddingProperty} We have seen that a hypersurface $M$ inside an
$\SU(3)$-holonomy 6-manifold $N$ is naturally endowed with a hypo structure $P$; hypo manifolds $(M,P)$ which occur
this way are said to have the embedding property. In this section we use Cartan-K\"ahler theory to prove that if $M$ is
real analytic and $P$ is a real analytic hypo structure, then $(M,P)$ has the embedding property.
 
 The analogous
problem can be studied in dimension 6, where one considers half-flat hypersurfaces inside $\Gtwo$-holonomy manifolds:
in \cite{Hitchin:StableForms}, Hitchin regarded half-flat structures as defined by closed forms $\sigma$ and $\rho$,
corresponding to $\omega^2/2$ and $\Re\Psi$ in our language. Having introduced a symplectic structure on the
infinite-dimensional affine space given by the product of the De Rham cohomology classes $[\sigma]$ and $[\rho]$,
Hitchin reduced the embedding property to a problem of existence of integral lines for a Hamiltonian vector field on
$[\sigma]\times[\rho]$.
This gives us an alternative approach to establish the embedding property, which we shall not discuss  in this paper. Nevertheless, we choose
to start this section by reformulating the embedding property in a more explicit way, in the spirit of
\cite{Hitchin:StableForms}:

\begin{proposition} A compact hypo manifold  $(M,\alpha,\omega_i)$ has the embedding property if and only if
$(\alpha,\omega_i)$
belongs to a one-parameter family of hypo structures  $(\alpha(t),\omega_i(t))$ satisfying the evolution equations
\begin{equation}\label{eqn:HypoEvolution}\left\{\begin{aligned}
\partial_t\omega_1&=-d\alpha\\
\partial_t(\omega_2\wedge\alpha)&=-d\omega_3\\
\partial_t(\omega_3\wedge\alpha)&=d\omega_2
\end{aligned}\right.\end{equation}
for all $t\in (a,b)$; then
\begin{equation*}\begin{aligned}
\omega&=\alpha\wedge dt+\omega_1\\
\Psi&=(\omega_2+i\omega_3)\wedge(\alpha+idt)
\end{aligned}\end{equation*}
defines an integrable $\SU(3)$-structure on $M\times(a,b)$.
\end{proposition}
\begin{proof}
The ``if'' part is obvious. Conversely, given a compact oriented hypersurface $M$ embedded in a $\SU(3)$-holonomy
6-manifold $N$, an embedding \[\phi\colon M\times(a,b)\ni (x,t)\to \exp_x(t \nu)\in N\] is defined for some interval
$(a,b)$, where $\nu\in T_xM$ is the unit normal compatible with the orientations. This gives a one-parameter family of
embeddings
\[\phi_t\colon M\ni x\to \phi(x,t)\in N\;,\] such that $\phi_0$ coincides with the original embedding.
The corresponding one-parameter family of hypo structures on $M$ evolves according to \eqref{eqn:HypoEvolution}.
\end{proof}
\begin{remark}
A hypo structure is Einstein-Sasaki if and only if the components \eqref{eqn:conicalSU3fromSU2} satisfy
\eqref{eqn:HypoEvolution}; similarly, nearly-K\"ahler half-flat structures are characterized by the evolution being
conical \cite{Hitchin:StableForms}. In this sense, evolution theory is a generalization of the construction of a
manifold with a parallel spinor as the cone on a manifold with a Killing spinor (see \cite{Bar}), where the cone is
replaced by more complicated evolution equations and the Killing spinor by a generalized Killing spinor.
\end{remark}

The rest of this section consists in the proof of the embedding theorem; for details on Cartan-K\"ahler theory we refer
to \cite{BryantEtAl}. Let $M$ be a real analytic five-manifold with real analytic forms $(\alpha,\omega_i)$ defining a
hypo structure. Define an embedding
\[\iota\colon M\ni x\to (x,0)\in N=M\times\R\;;\]
let $\pi\colon F\to N$ be the principal bundle of frames and $S=F/\SU(3)$. By Proposition \ref{prop:SU2fromSU3}, the
hypo structure on $M$ induces an $\SU(3)$-structure on $N$; by restriction, a section $f_S\in\Gamma(M, \iota^*S)$ is
defined. On a small open set $U\subset M$, $f_S$ lifts to a section $f$ of $\iota^*F$. The image $X_U=f(U)$ is a real
analytic
submanifold of $F$.

 Let $\mathcal{I}\subset\Omega(F)$ be the differential graded ideal generated by $d\omega$ and
$d\psi^\pm$, where
\begin{align*} \omega&=\eta^1\wedge\eta^2+\eta^3\wedge\eta^4+\eta^5\wedge \eta^6\;,\\
\psi^+ +i\psi^-&=(\eta^1+i\eta^2)\wedge(\eta^3+i\eta^4)\wedge(\eta^5+i\eta^6)\;,
\end{align*}
$\eta=(\eta^1,\dotsc,\eta^6)$ being the tautological form. This choice is consistent with \eqref{eqn:referenceformsSU3}, and
it differs from the choice of \cite{Bryant:Calibrated} by a permutation of indices. By construction, $X_U$ is an
integral manifold for $\mathcal{I}$, namely $\iota^*d\omega=0=\iota^*d\psi^\pm$ holds. Proving the embedding property
on $U$ amounts to finding an integral manifold of dimension 6 containing $X_U$ which is transverse to the fibres of
$F$; in order to obtain a global result, we shall
think of $\mathcal{I}$ as a differential system on $S$ and extend its integral manifold  $f_S(M)$.

Let $V_n(\mathcal{I},\pi)$ be the set of $n$-dimensional integral elements of $\mathcal{I}$ which do not
intersect $\ker\pi_*$. The following lemma states that $(\mathcal{I},\pi)$ is involutive, and is proved in
\cite{Bryant:Calibrated}:
\begin{lemma}[Bryant]
\label{lemma:Bryant} Every $E_6$ in $V_6(\mathcal{I},\pi)$ is the terminus of a regular flag
\[E_0\subset E_1\subset E_2\subset E_3\subset E_4\subset E_5\subset E_6\;,\]
where
\begin{align*}
E_5&=\{v\in E_6\mid \eta^6(v)=0\} & E_4&=\{v\in E_5\mid \eta^4(v)=0\}\\
 E_3&=\{v\in E_4\mid\eta^2(v)=0\} & E_2&=\{v\in E_3\mid\eta^5(v)=0\}\\
E_1&=\{v\in E_2\mid\eta^3(v)=0\} & E_0&=\{0\}\subset E_1
\end{align*}
\end{lemma}
Lemma \ref{lemma:Bryant} per se guarantees local existence of a Calabi-Yau structure on $N$, in the guise of a
6-dimensional integral manifold transverse to the fibres of $F$; in order to prove the embedding property locally, we
need such an integral manifold to contain $X_U$. In other words, we require the local integrable
$\SU(3)$\nobreakdash-structure to give back the starting hypo structure on $M=M\times\{0\}\subset N$ through
Proposition \ref{prop:SU2fromSU3}. To achieve this, we have to show that every $T_uX_U$ can be extended to a
6\nobreakdash-dimensional
$E_6$ satisfying the hypothesis of  Lemma \ref{lemma:Bryant}.

Recall that the polar space of an integral element $E$ is the union of all integral elements containing $E$. Define
$\gl(6,\R)_5$ (isomorphic to $\R^6$) to be the space of six by six matrices whose first five columns are zero.

\begin{lemma}
\label{lemma:PolarSpace} The tangent space to $X_U$ at a point $u$ is contained in an element $E_6$ of
$V_6(\mathcal{I},\pi)$. Moreover, the polar space of $T_uX_U$ is the direct sum of $E_6$ and a vertical vector space
which is mapped to
\[\gl(6,\R)_5\oplus \su(3)\subset \gl(6,\R)\] by any connection form. \end{lemma}

\begin{proof}
The hypo structure defines a  Riemannian metric on $M$; consider the product metric on $N$. Let $\sigma$ be the
connection form of the Levi-Civita connection on $F$; then $D\eta=0$ yields the structure equation
\[d\eta^i=\eta^j\wedge\sigma^i_j\;.\]
Now let $x\in M$, $u=f(x)$, $E=\im f_*\subset T_uF$. We must find a vector $v$ in $T_uF$ such that $E_6=E\oplus\langle
v\rangle$ is an integral element for $\mathcal{I}$ which is transverse to $\ker\pi_*$. Let $e_1,\dotsc,e_5$ be the
basis of $E$ characterized by $\eta_i(e_j)=\delta_{ij}$. By construction $\eta_6(E)=0$. Since we are allowed to modify
$v$ by a combination of the $e_i$, we can assume that $v$ satisfies
\[\eta_1(v)=\dotsb=\eta_5(v)=0\;;\] moreover by rescaling we can assume $\eta_6(v)=\frac{1}{2}$. Since we are using the product
metric, we have $\sigma^i_6=0=\sigma^6_i$ on $E$ for $i=1,\dotsc,6$.

We must solve
\begin{align*}
d\omega(e_i,e_j,v)&=0\;\forall\, i,j=1,\dotsc, 5\;, & d\psi^\pm(e_i,e_j,e_k,v)&=0\;\forall\, i,j,k=1,\dotsc, 5\;;
\end{align*}
these equations split up into eight independent blocks. Writing $\sigma^i_j$ for $\sigma^i_j(v)$, the first block is:

\begin{equation*}\left\{
\begin{aligned}
\sigma_1^1+\sigma_2^2 &=  d\eta^5(e_1,e_2) \\ \sigma_3^3+\sigma_4^4 &=   d\eta^5(e_3,e_4)\\
\sigma_2^2+\sigma_4^4+\sigma_5^5&=d\eta^1(e_2,e_5)+d\eta^3(e_4,e_5) \\
\sigma_1^1+\sigma_4^4+\sigma_5^5&=-d\eta^2(e_1,e_5)+d\eta^3(e_4,e_5)\\
\sigma_2^2+\sigma_3^3+\sigma_5^5&=d\eta^1(e_2,e_5)-d\eta^4(e_3,e_5) \\
\sigma_1^1+\sigma_3^3+\sigma_5^5&=-d\eta^2(e_1,e_5)-d\eta^4(e_3,e_5)
\end{aligned}\right.
\end{equation*}
Without even using  the fact that $E$ is an integral element, the above can be reduced to five independent, compatible
equations. Similarly, the following block can be reduced to three independent, compatible equations:
\begin{equation*}\left\{
\begin{aligned}
\sigma_1^2-\sigma_4^3+\sigma_5^6&=d\eta^1(e_1,e_5)+d\eta^4(e_4,e_5)\\
-\sigma_2^1+\sigma_3^4+\sigma_5^6&=d\eta^2(e_2,e_5)+d\eta^3(e_3,e_5)\\
\sigma_1^2+\sigma_3^4+\sigma_5^6&=d\eta^1(e_1,e_5)+d\eta^3(e_3,e_5)\\
\sigma_2^1+\sigma_4^3-\sigma_5^6&=-d\eta^2(e_2,e_5)-d\eta^4(e_4,e_5)
\end{aligned}\right.
\end{equation*}
Then we have four similar blocks:
\begin{align*}
&\left\{ \begin{aligned}
-\sigma_1^5+\sigma_2^6&=d\eta^1(e_1,e_2)-d\eta^4(e_2,e_4)+d\eta^3(e_1,e_4)\\
-\sigma_1^5+\sigma_2^6&=-d\eta^4(e_1,e_3)+d\eta^1(e_1,e_2)+d\eta^3(e_3,e_2)\\
\sigma_5^1+\sigma_2^6&=  d\eta^5(e_5,e_2)\\
\end{aligned}\right.\\
&\left\{ \begin{aligned}
-\sigma_3^5+\sigma_4^6&=d\eta^2(e_1,e_3)+d\eta^1(e_1,e_4)+d\eta^3(e_3,e_4)\\
-\sigma_3^5+\sigma_4^6&=d\eta^1(e_3,e_2)+d\eta^2(e_2,e_4)+d\eta^3(e_3,e_4)\\
\sigma_5^3+\sigma_4^6 &=  d\eta^5(e_5,e_4)\\
\end{aligned}\right.\\
&\left\{
\begin{aligned}
\sigma_2^5+\sigma_1^6&=-d\eta^4(e_3,e_2)-d\eta^2(e_1,e_2)-d\eta^3(e_1,e_3)\\
-\sigma_2^5-\sigma_1^6&=d\eta^2(e_1,e_2)+d\eta^4(e_1,e_4)+d\eta^3(e_2,e_4)\\
\sigma_5^2-\sigma_1^6&=d\eta^5(e_1,e_5)\\
\end{aligned}\right.\\
&\left\{
\begin{aligned}
-\sigma_4^5-\sigma_3^6&=d\eta^2(e_3,e_2)+d\eta^4(e_3,e_4)-d\eta^1(e_2,e_4)\\
\sigma_4^5+\sigma_3^6&=d\eta^1(e_1,e_3)-d\eta^4(e_3,e_4)-d\eta^2(e_1,e_4)\\
\sigma_5^4-\sigma_3^6 &=  d\eta^5(e_3,e_5)\\
\end{aligned}\right.
\end{align*}
that are compatible if and only if
\begin{equation}
\label{eqn:lemma:1} \left\{\begin{aligned}
d\eta^4(e_1,e_3)+d\eta^4(e_4,e_2)+d\eta^3(e_1,e_4)+d\eta^3(e_2,e_3)&=0\\
d\eta^2(e_1,e_3)+d\eta^2(e_4,e_2)+d\eta^1(e_1,e_4)+d\eta^1(e_2,e_3)&=0\\
d\eta^4(e_1,e_4)+d\eta^4(e_2,e_3)-d\eta^3(e_1,e_3)-d\eta^3(e_4,e_2)&=0\\
d\eta^2(e_1,e_4)+d\eta^2(e_2,e_3)-d\eta^1(e_1,e_3)-d\eta^1(e_4,e_2)&=0
\end{aligned}\right.
\end{equation}
Recall that $M$ is hypo, so the forms
\begin{align*}
\psi_2&=\eta^{135}+\eta^{425}\;, & \psi_3&=\eta^{145}+\eta^{235}
\end{align*}
are closed when pulled back to $M$ via $f$. We have
\begin{align*}
\eta_4\wedge\psi_2+\eta_3\wedge\psi_3&=0\;,& \eta_2\wedge\psi_2+\eta_1\wedge\psi_3&=0\;,\\
\eta_4\wedge\psi_3-\eta_3\wedge\psi_2&=0\;, & \eta_2\wedge\psi_3-\eta_1\wedge\psi_2&=0\;;
\end{align*}
pulling back to $M$ and taking $d$, we find that \eqref{eqn:lemma:1} is satisfied.\\
The remaining two blocks are
\begin{align*}
&\left\{
\begin{aligned}
-\sigma_3^2+\sigma_1^4 &= -d\eta^5(e_1,e_3) \\
\sigma_4^1-\sigma_2^3&=  -d\eta^5(e_2,e_4)\\
-\sigma_2^3-\sigma_1^4&=d\eta^4(e_2,e_5)-d\eta^3(e_1,e_5)\\
\sigma_4^1+\sigma_3^2&=d\eta^1(e_3,e_5)-d\eta^2(e_4,e_5)\\
\end{aligned}\right.\\
&\left\{
\begin{aligned}
\sigma_4^2+\sigma_1^3&=  d\eta^5(e_1,e_4)\\
\sigma_3^1+\sigma_2^4&=  -d\eta^5(e_2,e_3)\\
-\sigma_2^4+\sigma_1^3&=-d\eta^4(e_1,e_5)-d\eta^3(e_2,e_5)\\
-\sigma_3^1+\sigma_4^2&=d\eta^2(e_3,e_5)+d\eta^1(e_4,e_5)\\
\end{aligned}\right.
\end{align*} that are compatible if and only if
\begin{align*}
-d\eta^5(e_1,e_3)\!+\!d\eta^5(e_2,e_4)\!+\!d\eta^4(e_2,e_5)\!-\!d\eta^3(e_1,e_5)\!+\!d\eta^1(e_3,e_5)\!-\!d\eta^2(e_4,e_5)&=0\\
-d\eta^5(e_1,e_4)\!-\!
d\eta^5(e_2,e_3)\!-\!d\eta^4(e_1,e_5)\!-\!d\eta^3(e_2,e_5)\!+\!d\eta^2(e_3,e_5)\!+\!d\eta^1(e_4,e_5)&=0
\end{align*}
which can be rewritten as
\begin{align*}
f^*\left(-d\eta^5\wedge\psi_2 + d(\eta^{12}+\eta^{34})\wedge(\eta^{14}+\eta^{23})\right)&=0\\
f^*\left(-d\eta^5\wedge\psi_3 - d(\eta^{12}+\eta^{34})\wedge(\eta^{13}+\eta^{42})\right)&=0
\end{align*}
and follow from $\eta^5\wedge\psi_i=0=f^*d(\eta^{12}+\eta^{34})$.

Summing up, we have a system of 22 compatible, independent equations. This means that there exists some $v$
in the polar space of $E$ with $\eta^6(v)=\frac{1}{2}$. The characterization of this polar space in terms of
$E_6$ is obtained by repeating the calculations assuming $\eta_6(v)=0$; the resulting equations are obtained
from the ones above by setting the right-hand sides to zero.
\end{proof}
We can now prove the embedding theorem.
\begin{theorem}
\label{thm:EmbeddingProperty}  Let $M$ be a real analytic manifold with a real analytic hypo structure $P$. Then
$(M,P)$ has the embedding property.
\end{theorem}
\begin{proof}
Recall that locally, one can lift $f_S$ to a section $f$ with image $X_U$; let $u$ be a point of $X_U$. By Lemma
\ref{lemma:PolarSpace}, there is an element $E_6$ of $V_6(\mathcal{I},\pi)$ containing $T_uX$; by Lemma
\ref{lemma:Bryant}, $E_6$ is the terminus of a regular flag with $E_5=T_uX_U$. As a consequence, $T_uX_U$ is regular;
since this is true for all $u$, $X_U$ is regular. It follows that $f_S(U)$ is regular; since this is true for all $U$,
we can conclude that $f_S(M)$ is regular.
Now let $\SU(2)$ act on $\gl(6,\R)$ by conjugation. Observe that $\gl(6,\R)_5\oplus\su(3)$ is closed under $\SU(2)$
action; its orthogonal complement has dimension $22$, and will be denoted by $W_{22}$.

By construction $W_{22}$ does not intersect $\su(3)$; this implies that for a suitable neighbourhood $U_{22}$ of $0$ in
$W_{22}$ the map
\[U_{22}\times \SU(3)\ni (x,a)\to (\id+x)a\in \GL(6,\R)\] is an embedding. We can assume that $U_{22}$ is closed under
$\SU(2)$ action. Now consider the product $\SU(2)$-structure $P_N$ on $N$, induced by the hypo $\SU(2)$\nobreakdash-\hspace{0pt}structure on
$M$; recall that $P_N/\SU(3)\subset S$ contains $f_S(M)$. Now write
\[\mathcal{G}_0=P_N\times_{\SU(2)} (\id+U_{22})\subset P_N\times_{\SU(2)} \GL(6,\R)=\mathcal{G}\;;\]
identifying $\mathcal{G}$ with $F\times_{\GL(6,\R)}\GL(6,\R)$ we can define the natural bundle map
\begin{align*}
  \rho\colon F\times_N\mathcal{G}&\to F\\
 (u,[u,g])  &\to ug
\end{align*}
where $F\times_N\mathcal{G}$ is the fibre (i.e. fibrewise) product of $F$ and $\mathcal{G}$ over $N$. A bundle map $\rho\colon S\times_N\mathcal{G}\to S$ is induced. Let
\[Y=\rho\bigl((P_N/\SU(3))\times_N\mathcal{G}_0\bigr)\;;\] then $Y$ is a 28-dimensional real analytic submanifold of
$S$, which at each point $x$ of $f_S(M)$ is transverse to the polar space of $T_x(f_S(M))$. Since $f_S(M)$ is regular
and the codimension of $Y$ in $S$ is 6, coinciding with the extension rank of $f_S(M)$, we can apply the
Cartan-K\"ahler theorem and obtain a 6\nobreakdash-dimensional integral manifold for $\mathcal{I}$ on $S$ containing
$f_S(M)$. By construction, this 6-dimensional manifold is transverse to the fibres of $\pi$, and so it defines an
integrable $\SU(3)$\nobreakdash-structure on a neighbourhood of $M$ in $N$.
\end{proof}
\begin{remark} Even locally, Theorem \ref{thm:EmbeddingProperty} does not fully answer the embedding problem, because
  Calabi-Yau manifolds admit non-analytic hypersurfaces. What is
  worse, there exist non-analytic hypo structures on an open set of
  $\R^5$ which do not satisfy the embedding property.  Such examples
  can be constructed by ensuring that the intrinsic torsion has
  $g_2^3$ constant. This leads to a contradiction since the invariant
  $g_2^3$ represents the trace of the Weingarten tensor of a potential
  embedding, and a constant mean curvature hypersurface in a Calabi-Yau
  manifold is necessarily real analytic \cite{Bryant:Private}.
 \end{remark}

\section{Hypo nilmanifolds}
\label{sec:nilmanifolds} It is well known that nilmanifolds do not admit Einstein-Sasaki structures; in fact,
Einstein-Sasaki manifolds have finite fundamental group, and therefore $b_1=0$, which cannot occur for
nilmanifolds. Not surprisingly, most 5\nobreakdash-nilmanifolds do admit (invariant) hypo structures. Indeed,
consider
\[M=\Gamma\backslash G\;,\] where $G$ is a 5-dimensional nilpotent group, $\Gamma$ a discrete subgroup of $G$
and $M$ is compact; an invariant structure on the nilmanifold $M$ is a structure which pulls back to a
left-invariant structure on $G$. With this setting in mind, we define:

\begin{definition}
A hypo structure on  $\lie{g}$ is a quadruplet $(\alpha,\omega_i)$ of forms on $\lie{g}$ satisfying Proposition
\ref{prop:SU2forms} and Equation \ref{eqn:hypo}.
\end{definition}
\begin{remark}
Since Lie groups are real analytic, and invariant forms are real analytic, every hypo structure appearing in the
classification to follow satisfies the embedding property. \end{remark}
 Borrowing notation from \cite{Salamon:ComplexStructures}, we
represent Lie algebras using symbolic expressions such as $(0,0,0,0,12)$, which represents a Lie algebra with a basis
$e^1,\dotsc, e^5$ such that $de^i=0$ for $i=1,\dotsc,4$, and $de^5=e^{12}$.
\begin{theorem}
\label{thm:NilpotentHypo} The nilpotent 5-dimensional Lie algebras not admitting a hypo structure are $(0,0,12,13,23)$,
$(0,0,0,12,14)$ and $(0,0,12,13,14+23)$.
\end{theorem}
Using the classification of five-dimensional nilpotent Lie algebras, Theorem~\ref{thm:NilpotentHypo} is equivalent to
the following table:
\[\begin{array}{|r|r|r|r|}\hline
\lie{g}             &\text{step}    &b_2    &\text{Admits hypo}\\
\hline
0,0,12,13,14+23     &4              &3        &\text{no}\\
0,0,12,13,14        &4              &3        &\text{yes}\\
0,0,12,13,23        &3              &3        &\text{no}\\
0,0,0,12,14         &3              &4        &\text{no}\\
0,0,0,12,13+24      &3              &4        &\text{yes}\\
0,0,0,12,13         &2              &6        &\text{yes}\\
0,0,0,0,12+34       &2              &5        &\text{yes}\\
0,0,0,0,12          &2              &7        &\text{yes}\\
0,0,0,0,0           &1              &10       &\text{yes}\\
\hline
\end{array}\]
\begin{remark}
It follows from this table that among nilmanifolds $\Gamma\backslash G$ there are non-trivial examples of compact hypo manifolds with odd-dimensional $b_1$.
\end{remark}

We start with a list of examples of hypo structures on nilpotent Lie algebras which do admit such structures; we shall
then prove that any nilpotent Lie algebra with a hypo structure must be one of these. Theorem \ref{thm:NilpotentHypo}
will then follow from the classification of 5-dimensional nilpotent Lie algebras, which is not otherwise used.
\begin{itemize}
\item $(0,0,12,13,14)$ has a hypo structure given by
\begin{align*}\alpha&=e^1 & \omega_1&=e^{25}+e^{43} & \omega_2&=e^{24}+e^{35} & \omega_3&=e^{23}+e^{54}\end{align*}
\item $(0,0,0,12,13+24)$ has a one-parameter family of hypo structures given by
\begin{align*}\alpha&=e^1+e^5 & \omega_1&=e^4\wedge(-c\,e^2-e^3)+e^{25}\\
\omega_2&=e^{42}+e^5\wedge(-c\,e^2-e^3) & \omega_3&=e^{45}+(-c\,e^2-e^3)\wedge e^2\end{align*}
\item $(0,0,0,12,13)$ has a hypo structure given by
\begin{align*}\alpha&=e^1 & \omega_1&=e^{35}+e^{24} & \omega_2&=e^{32}+e^{45} & \omega_3&=e^{34}+e^{52}\end{align*}
 Taking the product of this nilmanifold with a circle, one obtains the half-flat symplectic structure in
\cite{Giovannini} (see Section~\ref{sec:examples}).
\item $(0,0,0,0,12+34)$ has hypo structures given by
\begin{align*}\alpha&=e^5 & \omega_1&=e^{12}+e^{34} & \omega_2&=e^{13}+e^{42} & \omega_3&=e^{14}+e^{23}\\
\alpha&=e^5 & \omega_1&=e^{12}-e^{34} & \omega_2&=e^{13}-e^{42} & \omega_3&=e^{14}-e^{23}\end{align*} These structures
arise as circle bundles over the hyperk\"ahler torus.
\item $(0,0,0,0,12)$ has hypo structures given by
\begin{align*}\alpha&=e^1 & \omega_1&=e^{25}+e^{34} & \omega_2&=e^{23}+e^{45} & \omega_3&=e^{24}+e^{53}\\
\alpha&=e^5 & \omega_1&=e^{12}+e^{34} & \omega_2&=e^{13}+e^{42} & \omega_3&=e^{14}+e^{23}\\
\alpha&=e^2-e^5 & \omega_1&=e^{34}+e^{15} & \omega_2&=e^{31}+e^{54} & \omega_3&=e^{35}+e^{41}\end{align*}
\item Every $\SU(2)$-structure on $(0,0,0,0,0)$ is hypo.
\end{itemize}

The rest of this section consists of the proof of the theorem. From now on, assume that $\lie{g}$ is a non-trivial
nilpotent Lie algebra carrying a hypo structure. Since $\lie{g}$ is nilpotent, one can fix a filtration of vector
spaces $V^i$, $\dim V^i=i$, such that \[V^1\subset V^2\subset \dotsb\subset V^5=\lie{g^*}, \quad
d(V^i)\subset\Lambda^2V^{i-1}\;.\] This filtration can be chosen so that $V^i=\ker d$ for some $i$; in particular, one
has $V^2\subset\ker d\subset V^4$. Note that the first Betti number $b_1$ is the dimension of $\ker d$.

It is convenient to distinguish three cases, according to whether $\alpha$ lies in $V^4$, $(V^4)^\perp$ or neither.

\subsection{First case}
We first consider the case when $\alpha$ is in $V^4$.
\begin{theorem}
If $\alpha$ is in $V^4$, then $\lie{g}$ is either $(0,0,0,0,12)$, $(0,0,0,12,13)$, or $(0,0,12,13,14)$.
\end{theorem}
\begin{proof}
Fix a unit $e^5$ in $(V^4)^\perp$ and apply Corollary \ref{cor:referenceformsSU2} to obtain a coframe $e^1,\dotsc,e^5$
such that
\begin{align*}
\alpha&=e^1\;, & \omega_1&=e^{25}+e^{34}\;,& \omega_2&=e^{23}+e^{45}\;, & \omega_3&=e^{24}+e^{53}\;.
\end{align*}
From $d\omega_1=0$, it follows that
 \begin{equation}\label{eqn:thm:b1lessthan4:1} e^2\wedge de^5-de^{34}=de^2\wedge e^5\;;
 \end{equation} since the left-hand side lies in $\Lambda^3 V^4$ and the right-hand side lies in $e^5\wedge\Lambda^2
V^4$, both must vanish and $de^2=0$. Using the fact that $\alpha\wedge\omega_2$ and $\alpha\wedge\omega_3$ are closed,
we obtain:
 \begin{align*}
0&=de^{123}+de^{14}\wedge e^5+e^{14}\wedge de^5\\
0&=de^{124}-de^{13}\wedge e^5-e^{13}\wedge de^5
\end{align*}
Therefore $e^{13}$ and $e^{14}$ are closed; since $e^2$ is closed as well, we get
\begin{align}
\label{eqn:thm:b1lessthan4:2} 0&=e^{14}\wedge de^5 \;,& 0&=e^{13}\wedge de^5\;.
\end{align}

Suppose first that $e^1$ is not in $V^3$. By dimension count in $V^4$, $\langle e^3,e^4\rangle$ has non-zero
intersection with $V^3$, and we can therefore rotate $e^3$ and $e^4$ to get a different hypo structure with $e^3\in
V^3$; then $de^{13}=0$ implies that $e^3$ is closed and $e^3\wedge de^1=0$. Similarly, $de^4$ is a multiple of $de^1$:
write $e'=a\,e^1+b\,e^4\in V^3$, where $b$ must be non-zero; then \[d(e^1\wedge e')=b\,de^{14}=0\;,\] so $e'$ is
closed. It follows that $e^{34}$ is closed; then \eqref{eqn:thm:b1lessthan4:1} becomes
\[e^2\wedge de^5=0\;,\] showing that, because of \eqref{eqn:thm:b1lessthan4:2}, up to a constant $de^5$ equals
$e^{12}$, which is therefore closed. Thus, all 2-forms on $V^4$ are closed, implying that $V^4$ is trivial as a Lie
algebra; consequently, $\lie{g}=(0,0,0,0,12)$.

Assume now that $e^1$ is in $V^3$; we can rotate $e^3$ and $e^4$ to get $e^3$ in $V^3$, $e^4$ in $(V^3)^\perp$. From
$de^{14}=0$ we find that $e^1$ is closed and $de^{4}\wedge e^1=0$. Wedging the left-hand side of
\eqref{eqn:thm:b1lessthan4:1} with $e^1$ and using $de^{13}=0$, we see that
\[e^{12}\wedge de^5=e^1\wedge de^{34}=e^{14}\wedge de^3-e^{13}\wedge de^4=e^{34}\wedge de^1=0\;.\] Together with
\eqref{eqn:thm:b1lessthan4:2}, this implies that $de^5$ is in $\langle e^{12},e^{13},e^{14}\rangle$. Now consider the
endomorphism $f$ of $\alpha^\perp$ defined by $e^1\wedge f(\eta)=d\eta$; its matrix with respect to
$\{e^2,e^3,e^4,e^5\}$ is strictly upper triangular. Its Jordan canonical form is therefore one of
\[
\begin{pmatrix}
0&0&0&0\\
0&0&0&0\\
0&0&0&1\\
0&0&0&0
\end{pmatrix}\;,
\begin{pmatrix}
0&1&0&0\\
0&0&0&0\\
0&0&0&1\\
0&0&0&0
\end{pmatrix}\;,
\begin{pmatrix}
0&1&0&0\\
0&0&1&0\\
0&0&0&1\\
0&0&0&0
\end{pmatrix}\; \text{ or }
\begin{pmatrix}
0&0&0&0\\
0&0&1&0\\
0&0&0&1\\
0&0&0&0
\end{pmatrix}\;;
\]
the first three cases give $\lie{g}=(0,0,0,0,12)$, $(0,0,0,12,13)$, $(0,0,12,13,14)$ respectively, and we must show
that the last case cannot occur. Indeed, in the last case $\lie{g}=(0,0,0,12,14)$; then
\[\ker d=V^3=\langle e^1,e^2,e^3\rangle\;,\] and on the other hand
\begin{align*}
e^2\wedge de^5&\notin \Lambda^3V^3\;, & de^{34}=e^3\wedge de^4\in\Lambda^3V^3\;,
\end{align*}
contradicting \eqref{eqn:thm:b1lessthan4:1}.
\end{proof}

\subsection{Second case}
The key tool to classify the remaining hypo structures is the following lemma, which shows that $\alpha$ is orthogonal
to $V^4$ if and only if $\omega_2$, $\omega_3$ are closed and $b_1=4$.

\begin{lemma}
\label{lemma:alphaorthogonal} Let $\alpha\notin V^4$. Then
 \begin{enumerate}
 \item If all $\omega_i$ are closed, $\alpha$ is orthogonal to $V^4$.
 \item If $\alpha$ is orthogonal to $V^4$, $V^4=\ker d$. In particular, all $\omega_i$ are closed.
 \end{enumerate}
\end{lemma}
\begin{proof}
Let $e^5$ be a unit form in $(V^4)^\perp$. By hypothesis, $\alpha=a\,e^5+\eta$ where $\eta$ is in $V^4$ and $a\neq 0$.
To prove the first statement, suppose $\eta$ is non-zero and let $e^4$ be a unit form in $\langle e^5,\eta\rangle$
orthogonal to $\alpha$. Using Corollary \ref{cor:referenceformsSU2} we can write
\begin{align*}
\omega_1&=e^{12}+e^{34}\;, & \omega_2&=e^{13}+e^{42}\;, & \omega_3&=e^{14}+e^{23}\;.
\end{align*}
The space $\langle e^1,e^2,e^3\rangle$ is orthogonal to $e^5$, and is therefore contained in $V^4$, whereas $e^4$ is
not. Since $\omega_1$ is closed,
\begin{equation}
\label{eqn:lemma:alphaorthogonal} e^3\wedge de^4-de^{12}=de^3\wedge e^4\;;
\end{equation}
both sides must then vanish, and so $e^3$ is closed; applying the same argument to $\omega_2$ and $\omega_3$ (which are
closed by hypothesis) one finds that $e^1$ and $e^2$ are also closed. From \eqref{eqn:lemma:alphaorthogonal} and its
analogues obtained using $\omega_2$ and $\omega_3$, it
follows that $de^4\wedge \langle e^1,e^2,e^3\rangle$ is trivial. Hence $e^4$ is closed and is therefore in $V^4$, which is absurd.

To prove the second assertion, let $\eta=0$, i.e. $\alpha=e^5$ (up to  sign). From \eqref{eqn:hypo}, it follows that
$\omega_2$ and $\omega_3$ are closed. Pick a unit $e^4$ in $(V^3)^\perp\cap V^4$, and define $e^1,e^2,e^3$ so as to
obtain \eqref{eqn:referenceformsSU2}; then $\langle e^1,e^2,e^3\rangle=V^3$. The same argument as above gives $V^4=\ker
d$.
\end{proof}
It is now easy to prove:
\begin{theorem}
If $\alpha$ is orthogonal to $V^4$, then $\lie{g}$ is one of \[(0,0,0,0,12)\;, \quad(0,0,0,0,12+34)\;.\]
\end{theorem}
\begin{proof}
By Lemma \ref{lemma:alphaorthogonal}, $V^4=\ker d$. Then either $d\alpha$ is simple (i.e. of the form $e\wedge f$), and one can choose a basis
$e^1,\dotsc,e^4$ of $V^4$ such that $d\alpha=e^{12}$, or it is not simple, and one can choose a basis such that
$d\alpha=e^{12}+e^{34}$.
\end{proof}

\subsection{Third case}
The last case is the one with $\alpha$ neither in $V^4$ nor in $(V^4)^\perp$. Lemma \ref{lemma:alphaorthogonal}
suggests that the span of $d\omega_2$ and $d\omega_3$ is relevant to the classification of hypo structures; we shall
use the dimension of $\langle d\omega_2,d\omega_3\rangle\cap\Lambda^3V^4$ to distinguish two subcases. In fact, we
shall prove that this dimension can only be $1$ or $2$.

\begin{theorem}
\label{thm:domegainV4} If $\alpha$ is neither in $V^4$ nor in $(V^4)^\perp$ and $d\omega_2$, $d\omega_3$ are in
$\Lambda^3V^4$, then $\lie{g}=(0,0,0,0,12)$.
\end{theorem}
\begin{proof}
Let $\alpha+\gamma$ be a generator of $(V^4)^\perp$ with $\gamma$ in $\alpha^\perp$, and let $k$ be the norm of
$\gamma$; then $\alpha-k^{-2}\gamma$ lies in $V^4$. Consider now the hypo structure obtained multiplying $\alpha$ by
$k$. Let $e^4=-k^{-1}\gamma$ and define
 \begin{align*}
 \eta&=\frac{1}{2}\left(\alpha+e^4\right)\;, &  \xi&=\frac{1}{2}\left(\alpha-e^4\right)\;;
\end{align*}
then $\xi$ generates $(V^4)^\perp$ and $\eta$ is in $V^4$. 
Using Corollary \ref{cor:referenceformsSU2} we can write
\begin{align*}
 \omega_1&=e^{12}+e^{34}\;, & \omega_2&=e^{13}+e^{42}\;, & \omega_3&=e^{14}+e^{23}\;.
\end{align*}
The space $\langle e^1,e^2,e^3\rangle$, being orthogonal to both $e^4$ and $\alpha$, is orthogonal to $\xi$, and is
therefore contained in $V^4$; on the other hand $e^4$ is not in $V^4$. Since $\omega_1$ is closed,
 \begin{align}
 \label{eqn:inV4:domega1}
 e^3\wedge de^4-de^{12}=de^3\wedge e^4\;;
 \end{align}
both sides must then vanish, and so $e^3$ is closed. Similarly, write
\begin{equation}
\label{eqn:inV4:domega23}
 \begin{aligned}
 d\omega_2&=de^{13}+de^4\wedge e^2-\eta\wedge de^2+\xi\wedge de^2\\
 d\omega_3&=de^{23}-de^4\wedge e^1+\eta\wedge de^1-\xi\wedge de^1
 \end{aligned}
\end{equation}
 So far we have only used the fact that $\alpha$ lies neither in $V^4$ nor in $(V^4)^\perp$. Writing \[\Lambda^3\lie{g}^*=\Lambda^3 V^4\oplus \xi\wedge\Lambda^2 V^4\;,\]
 the hypotheses on $d\omega_2$, $d\omega_3$ imply that $e^2$ and $e^1$ are closed.
By \eqref{eqn:hypo},
\begin{align}
\label{eqn:inV4:domega2alpha} 0&=d(\omega_2\wedge\alpha)=d(e^{13}\wedge (\eta+\xi)-2\,e^{2}\wedge\eta\wedge\xi)\\
\label{eqn:inV4:domega3alpha} 0&=d(\omega_3\wedge\alpha)=d(e^{23}\wedge (\eta+\xi)+2\,e^1\wedge\eta\wedge\xi)
\end{align}
Relative to $\Lambda^4\lie{g}^*=\xi\wedge\Lambda^3 V^4\oplus \Lambda^4V^4$, the first components give $e^1\wedge
d\eta=0$, $e^2\wedge d\eta=0$; using this, the second components give
\begin{align*}
(e^{13}-2\,e^2\wedge\eta)\wedge d\xi=0\\
(e^{23}+2\,e^1\wedge\eta)\wedge d\xi=0
\end{align*}
Since $e^1$, $e^2$, $e^3$  are closed, \eqref{eqn:inV4:domega1} becomes
 \begin{equation}
 \label{eqn:lemma:domegainV4:1}
 e^3\wedge d\eta=e^3\wedge d\xi\;.
 \end{equation} Wedging by $e^1$, $e^2$ we see that
$d\xi\wedge e^{13}$ and $d\xi\wedge e^{23}$ are zero, so our equations reduce to
 \begin{align*}
 (e^2\wedge\eta)\wedge d\xi&=0\;, & (e^1\wedge\eta)\wedge d\xi&=0\;.
 \end{align*}
Therefore, $d\xi$ lies in $\langle e^{12},e^3\wedge\eta\rangle$. Suppose that $\eta$ is not closed; then $d\eta$ is a
non-zero multiple of $e^{12}$, so $d^2=0$ implies that $d\xi$ must also be a multiple of $e^{12}$. Then
\eqref{eqn:lemma:domegainV4:1} implies that $e^4$ is closed, which
is absurd because $e^4$ cannot be in $V^4$.
Thus, $\eta$ is closed and from \eqref{eqn:lemma:domegainV4:1} we must have $d\xi=e^3\wedge\eta$ up to a constant,
which we can take to be $1$ by introducing a global scale factor. The choice of basis $\{e^3,\eta,e^1,e^2,e^4\}$ reveals
$\lie{g}$ to be $(0,0,0,0,12)$.
\end{proof}

Suppose now that the dimension of $\langle d\omega_2,d\omega_3\rangle\cap\Lambda^3V^4$ is one; then up to
rotating $\omega_2$ and $\omega_3$ we can assume that, say, $d\omega_2$ is in $\Lambda^3 V^4$. The following
lemma shows that this can always be done.
\begin{lemma}
\label{lemma:assumedomega2inV4} If $\alpha$ is neither in $V^4$ nor in $(V^4)^\perp$, up to rotating $\omega_2$ and
$\omega_3$ we can always assume that $d\omega_2$ is in $\Lambda^3V^4$.
\end{lemma}
\begin{proof}
We can proceed as in the proof of Theorem \ref{thm:domegainV4}. From \eqref{eqn:inV4:domega23} we
see that if $de^1$ and $de^2$ are linearly dependent, the statement holds; assume that they are independent.
Consider the symmetric bilinear form $B$ on $\Lambda^2V^4$ defined by
\[B(\alpha,\beta)e^{123}\wedge\eta=\alpha\wedge\beta\;;\] its signature is $(+,+,+,-,-,-)$.  By the classification of  nilpotent
 Lie algebras of dimension four, $V^4=(0,0,12,13)$; an explicit computation then shows that the space $Z_2$ of
closed 2-forms has dimension 4 and the signature of $B$ on $Z_2$ is $(0,0,+,-)$.
On the other hand, the components of \eqref{eqn:inV4:domega2alpha}, \eqref{eqn:inV4:domega3alpha} containing
$\xi$ give
\begin{align*}
d(e^{13}+2\,\eta\wedge e^2)&=0\;, & d(2\,e^1\wedge\eta+e^{23})&=0\;,
\end{align*}
giving a two-dimensional subspace in $Z_2$ on which $B$ is positive definite, which is absurd.
\end{proof}

\begin{theorem}
\label{thm:domega3notinV4} Let $\alpha\notin V^4$. If $d\omega_3$ is not in $\Lambda^3V^4$, then
\[\lie{g}=(0,0,0,12,13+24)\;.\]
\end{theorem}
\begin{proof}
Since $d\omega_3$ is not in $\Lambda^3V^4$, $\alpha$ cannot be orthogonal to $V^4$, because otherwise
$\omega_3$ would be in $\Lambda^2 V^4$. We can then proceed as in the proof of Theorem \ref{thm:domegainV4};
in fact, everything applies verbatim until the conclusion that $e^1$ and $e^2$ are closed, as in the present
case $e^2$ is closed but $e^1$ is not. Equations \ref{eqn:inV4:domega2alpha} and \ref{eqn:inV4:domega3alpha}
also hold; the vanishing of the $\xi\wedge\Lambda^3V^4$ component of \eqref{eqn:inV4:domega3alpha} shows that
$d(e^1\wedge\eta)$=0. This implies that $d\eta=k\,de^1$ for some constant $k$, for if one were to rotate
$\eta$ and $e^1$ in order to have $\eta\in V^3$, then $\eta$ would become closed. Rewrite
\eqref{eqn:inV4:domega1}, \eqref{eqn:inV4:domega2alpha} and \eqref{eqn:inV4:domega3alpha} as
\begin{align}
\label{eqn:domega3notinV4:1}(k\,e^3-e^2)\wedge de^1-e^3\wedge d\xi=0\\
\label{eqn:domega3notinV4:2}de^1\wedge(2\,k\,e^2+e^3)=0\\
\label{eqn:domega3notinV4:3}(e^{13}-2\,e^2\wedge\eta)\wedge d\xi=0\\
\label{eqn:domega3notinV4:4}k\,e^{23}\wedge de^1+(e^{23}+2\,e^1\wedge\eta)\wedge d\xi=0
\end{align}
Wedging \eqref{eqn:domega3notinV4:1} with $e^3$ shows that $de^1\wedge e^{23}$ is zero; wedging with $e^2$ shows that
$d\xi\wedge e^{23}$ is zero. Equation \ref{eqn:domega3notinV4:4} is therefore equivalent to the vanishing of
$d\xi\wedge e^1\wedge\eta$. From $d(e^1\wedge\eta)=0$ and \eqref{eqn:domega3notinV4:2}, we find that up to a non-zero
multiple
\[de^1=(2\,k\,e^2+e^3)\wedge(\eta-k\,e^1)\;.\] Then from \eqref{eqn:domega3notinV4:1} we get
 \begin{equation}
 \label{eqn:domega3notinV4:dxi}
d\xi\in \langle e^2\wedge(\eta-k\,e^1)\rangle\oplus e^3\wedge V^4\;;
\end{equation}
 write $d\xi=\sigma_1+\sigma_2$ accordingly and note that $\sigma_1$
cannot be zero, as in that case \eqref{eqn:domega3notinV4:1} would imply $(2k^2+1)\,e^{23}=0$. From
\eqref{eqn:domega3notinV4:dxi} and \eqref{eqn:domega3notinV4:4}, we know that $\sigma_2$ must be in $e^3\wedge\langle
e^1,\eta\rangle$; let
\[\sigma_2=c_1\,e^3\wedge(\eta-k\,e^1)+c_2\,e^{31}\;.\]
If $c_2$ is zero, \eqref{eqn:domega3notinV4:3} implies that $d\xi=de^1$ up to a multiple. So $\langle e^1,\xi\rangle$
has non-zero intersection with $\ker d\subset V^4$, contradicting the fact that $\langle e^1,\xi\rangle$ intersects
$V^4$
in $\langle e^1\rangle$ and $de^1\neq 0$.
Hence $c_2$ is not zero, and since $\sigma_1$ is closed, $\sigma_2$ must be closed as well; therefore $de^1\wedge
e^3=0$, implying $k=0$. We can then rescale everything so that, using \eqref{eqn:domega3notinV4:1} and
\eqref{eqn:domega3notinV4:3},
\begin{align*}
  d\eta&=0\;, & de^1&=2\,e^3\wedge\eta\;, & d\xi&=2\,e^2\eta+e^{13}+c_1\,e^3\wedge\eta\;.
 \end{align*}
Relative to the basis $\{-2\;\eta, e^3, -e^2-\frac{c_1}{2}\,e^3, e^1,e^4\}$ we see that
\begin{equation*}\lie{g}=(0,0,0,12,13+24)\;.\qedhere
\end{equation*}
\end{proof}

\section{Examples and applications}\label{sec:examples}
 In this final section, we pick out two of the nilmanifolds from the previous section, and investigate the resulting
geometrical structures using techniques developed above.

We first give an example of a hypo structure such that $A$
does not satisfy the Codazzi equation, proving that Theorem \ref{thm:EmbeddingProperty} is in some sense more general
than the theorem of B\"ar, Gauduchon and Moroianu \cite{BarGauduchonMoroianu}. It is based on the Lie algebra
\begin{align}\label{eqn:g1234}
\lie{g}=(0,0,0,12,13+24),
\end{align}
 though the technique is likely to extend to other cases.

\begin{example} Given \eqref{eqn:g1234}, consider the hypo structure
\begin{align*}
\alpha&=e^1+e^5\;, & \omega_1&=-e^{43}+e^{25} \;,& \omega_2&=e^{42}-e^{53} \;,& \omega_3&=e^{45}-e^{32}\;.
\end{align*}
In order to check that the Codazzi equation is not satisfied, we now reinterpret $A$ as a map $\Theta\colon P\to
\Sym(\R^5)$, which under \eqref{eqn:GeneralizedKillingTorsion} can be identified with the intrinsic torsion. For
convenience, we work on a left-invariant reduction of $P$ to the trivial group $\{e\}$, namely with the global frame
\begin{equation*}
e^4\;,\quad -e^3\;,\quad e^2\;,\quad e^5\;,\quad e^1+e^5\;; \end{equation*}
 then the connection form can be written as \[
\begin {pmatrix}0&0&\eta^{{4}}&0&\eta^{{3}}
\\\noalign{\medskip}0&0&0&\eta^{{4}}&-\eta^{{5}}\\\noalign{\medskip}-\eta^{{4}}
&0&0&-\eta^{{1}}&0\\\noalign{\medskip}0&-\eta^{{4}}&\eta^{{1}}&0&-\eta^{{2}}
\\\noalign{\medskip}-\eta^{{3}}&\eta^{{5}}&0&\eta^{{2}}&0\end{pmatrix}\;,\]
where $\eta^k$ is the $k$-th component of the tautological form. From \eqref{eqn:HypoConnectionForm}, the restriction of
$\Theta$ to the $\{e\}$-structure is
\[
\Theta=\begin{pmatrix}
0       &   0   &   0   &   0   &   -1\\
0       &   0   &   -1  &   0   &   0\\
0       &   -1  &   0  &   0   &   0\\
0       &   0   &   0   &   0   &   0\\
-1      &   0   &   0   &   0   &   0\\
 \end{pmatrix}
\]
and we can identify $\nabla\Theta$ with the $\Sym(\R^5)$-valued one-form
 \[
 \begin{pmatrix}
-2\,\eta^3      &   \eta^5-\eta^4   &   0       &   \eta^2    &    0\\
\eta^5-\eta^4 & 0                   &   0       &   -\eta^1   &   0\\
0               &   0                   &   0       &   \eta^4    & \eta^4-\eta^5\\
\eta^2       &   -\eta^1           &   \eta^4 &   0           & 0\\
0               &   0                   &   \eta^4-\eta^5 & 0   & 2\,\eta^3
 \end{pmatrix}\;.\]
The Codazzi equation is satisfied if and only if $\nabla\Theta$ is totally symmetric as a section of  $T^*M\otimes
T^*M\otimes T^*M$, which is clearly false.
\end{example}

\vskip5pt

The remainder of our study is based exclusively on the Lie algebra
\begin{align}\label{eqn:g1213}
\lie{g}=(0,0,0,12,13)
\end{align}
that also figures in the table of Theorem~\ref{thm:NilpotentHypo}.

\begin{example} Consider the hypo structure already defined on \eqref{eqn:g1213}, namely the one with
\begin{align}\label{eqn:hypo1g}
\alpha&=e^1\;, & \omega_1&=e^{35}+e^{24}\;, & \omega_2&=e^{32}+e^{45}\;, & \omega_3&=e^{34}+e^{52}.
\end{align}
To solve the evolution equations \eqref{eqn:HypoEvolution}, we look (with hindsight) for an orthonormal basis on
$M\times(a,b)$ of the form
\begin{equation}\label{eqn:basis}
\begin{aligned}
E^1&=f\,e^1, & E^2&=g\,e^2, & E^3&=g\,e^3,\\
E^4&=g^{-1}e^4, & E^5&=g^{-1}e^5, & E^6&=dt,
\end{aligned}
\end{equation}
where $f=f(t)$ and $g=g(t)$. More precisely, we deform the hypo structure so that \begin{align*}
\alpha(t)&=f(t)e^1\;,&\omega_1(t)&=e^{35}+e^{24}\;,\\\omega_2(t)&=g(t)^2e^{32}+g(t)^{\!-2}e^{45}\;,&
\omega_3(t)&=e^{34}+e^{52}.\end{align*}

The evolution equations \eqref{eqn:HypoEvolution} then boil down to \begin{align*} \partial_t (fg^2)&=-2\;,&
\partial_t(fg^{-2})&=0\;, & \partial_t f&=-g^{-2}\;.& \end{align*} Without regard to initial conditions, the second
equation implies that $f=cg^2$ with $c$ constant. The remaining equations both become $2cg^3\partial_tg=-1$, giving
\begin{align}\label{eqn:solutiong} g(t)=(p\kern2pt t+q)^{1/4}, \end{align} where $p$, $q$ are constants. Any values of
the latter with $p\neq 0$ give rise to a solution of the evolution equations, but for the case in hand $f(0)=1=g(0)$.
This forces $p=-2$, $q=1$, and we may take $a=-\infty$, $b=\frac12$.  One can check that the resulting metric
\[(1\!-\!2t)\ee1+(1\!-\!2t)^{\frac12}(\ee2\!+\!\ee3) +(1\!-\!2t)^{\!-\frac12}(\ee4 \!+\! \ee5) + dt\!\otimes\kern-1pt dt\] is
Ricci-flat. Further curvature calculations confirm that its holonomy group is indeed \textit{equal} to $\SU(3)$.
\end{example}

Let us consider further the geometry induced on $N=M\times(-\infty,\frac12)$ by the above solution. We can extend the
structure group $\SU(3)$ of $N$ to $\LieG{U}(3)$ by considering the natural action of $S^1$ on the orthonormal basis
\eqref{eqn:basis} induced by the complex structure.  Thus, $u=e^{i\theta}$ acts by simultaneous rotation in the planes
$\left<E^3,E^5\right>,\left<E^2,E^4\right>,\left<E^1,E^6\right>$; for example,
\begin{equation}\label{eqn:rotation}
\left\{\begin{array}{rcl}
u\cdot E^3 &=& E^3\cos\theta - E^5\sin\theta,\\[3pt]
u\cdot E^2 &=& E^2\cos\theta - E^4\sin\theta,\\[3pt]
u\cdot E^1 &=& E^1\cos\theta - E^6\sin\theta.
\end{array}\right.\end{equation}
This action extends to differential forms by setting \begin{equation}\label{eqn:rotforms} u\cdot E^{jk\cdots}=(u\cdot
E^j)\wedge (u\cdot E^k)\wedge\cdots \end{equation} By construction, $u$ leaves invariant the K\"ahler form
$\omega=E^{35}+E^{24}+E^{16}$, but acts by multiplication by $e^{3i\theta}$ on the holomorphic 3-form
\begin{equation}\label{eqn:Psi}
\Psi=(E^3+iE^5)\wedge(E^2+iE^4)\wedge(E^1+iE^6),
\end{equation}
which is consistent with \eqref{eqn:referenceformsSU3} given the choice of indices in \eqref{eqn:hypo1g}.

To give an explicit example of the induced action on 3-forms, let $\varepsilon=e^{2\pi i/3}$. Then
\begin{equation}\label{eqn:Eisenstein}\begin{array}{rcl}
\varepsilon\cdot 8E^{326} &=& -(E^3+\sqrt3\,E^5)\wedge(E^2+\sqrt3\,E^4)\wedge(E^6-\sqrt3\,E^1)\\[5pt] &=&
g^2 e^{23}\!\wedge dt - \sqrt3\,\alpha - 3\,\beta - 3\sqrt3\,e^{145}\;,
\end{array}\end{equation}
where
\begin{equation}\label{eqn:alphabeta} \begin{array}{rclcl}
\alpha &=& -E^{321}+E^{346}+E^{526} &=& \!g^4 e^{123}+e^{34}\!\wedge dt-e^{25}\!\wedge dt\\[5pt]
\beta  &=& -E^{341}-E^{521}+E^{546} &=& -g^2e^{134}+g^2e^{125}-g^{-2}e^{45}\!\wedge dt
\end{array}\end{equation}
One can associate to a simple 3-form $\gamma$ the 3-dimensional subspace
\[\gamma\ann=\{v\in\lie{g}:v\,\hook\,\gamma=0\}\] that it annihilates; for example \begin{equation}\label{326o}
V=\left<E_5,E_4,E_1\right>=(E^{326})\ann,\end{equation} where $(E_i)$ is the dual (orthonormal) basis of tangent
vectors. In this language, \eqref{eqn:rotforms} induces the natural action on subspaces.

It was shown by Giovannini \cite{Giovannini} that the nilmanifold $M\times S^1$ based on the 6-dimensional nilpotent Lie
algebra $(0,0,0,0,12,13)$ admits a \textit{tri-Lagrangian structure}, meaning that it is a symplectic manifold through
each point of which pass three mutually transverse Lagrangian submanifolds. The notation \eqref{eqn:alphabeta} is in
fact taken from \cite{Giovannini}. It leads to an even richer structure in the Calabi-Yau setting, based on

\begin{lemma}
 \label{eqn:simple3forms} For each fixed $u\in S^1$, the simple 3-form $u\cdot E^{326}$ is closed.
\end{lemma}

\begin{proof} Referring to \eqref{eqn:Psi} and \eqref{eqn:Eisenstein}, we see that the $S^1$ orbit containing $E^{326}$
spans the four-dimensional vector space
\[
\left<E^{326},\>E^{541},\>\Re\Psi,\>\Im\Psi\right>=\left<E^{326},\>\alpha,\>\beta,\>E^{541}\right>.
\]
But both $E^{326}=g^2e^{32}\wedge dt$ and $E^{541}=e^{541}$ are closed, as is $\Psi$. Alternatively, one can use
\eqref{eqn:solutiong} and \eqref{eqn:alphabeta}, still assuming $p=-2$, $q=1$,  to verify directly that
$d\alpha=0=d\beta$.\end{proof}

Observe from \eqref{eqn:Psi} that $E^{326}\wedge\Im\Psi=0$. This is equivalent to asserting that the restriction of
$\Im\Psi$ to the subspace \eqref{326o} is identically zero, or that $V$ is \textit{special Lagrangian}
\cite{CalibratedGeometries}. Setting $u=e^{i\theta}$, it also follows that
\begin{align*}
        (u\cdot E^{326})\wedge\Im\!\left[e^{3i\theta}\Psi\right]=0.
\end{align*}
Hence, $u\cdot V$ is special Lagrangian \textit{with phase $e^{-3i\theta}$} in the sense of
\cite{Joyce:SpecialLagrangian}. In particular, the subspaces $V,\,\varepsilon\cdot V,\,\varepsilon^2 V$, corresponding
to $\theta=n\pi/3$ with $n\in\mathbb{Z}$, are all special Lagrangian with the same phase. The geometrical structure
defined by such a triple of Lagrangian subspaces has special interest, since the group $\LieG{Sp}(n,\R)$ of linear
symplectic transformations on $\R^{2n}$ acts almost transitively on triples of mutually transverse Lagrangian
subspaces, with stabilizer $\LieG{O}(n)$. Indeed, the only invariant of an $\LieG{Sp}(n,\R)$-orbit is the Maslov index
or signature (of the $n\times n$ symmetric matrix expressing the third subspace relative to the first two in standard
form).

The invariant structure defined by Lemma~\ref{eqn:simple3forms} extends to each point of $N=M\times(-\infty,\frac12)$,
and closedness of a 3-form translates into integrability of the associated distribution. Thus, we can assert the
existence of an $S^1$ `pencil' of Lagrangian manifolds through each point of $N$. For each fixed phase, exactly three
of them will be special Lagrangian. Similar configurations of submanifolds, in the guise of `D6 branes at $\SU(3)$
angles', occur in various contexts in M-theory, in particular from the fixed points of a circle action on the
$G_2$-holonomy cone over $\SU(3)/T^2$ \cite{AtiyahWitten, AcharyaDenefHofmanLambert}.

\begin{example} We shall now exhibit compatibility of the same structure with a very different metric with holonomy
$\SU(3)$ defined nonetheless on the same nilmanifold. Given \eqref{eqn:g1213} again, consider the hypo structure for
which
\begin{align}\label{eqn:hypo2g}
\alpha&=e^1\;, & \omega_1&=e^{35}-e^{24}\;, & \omega_2&=-e^{32}+e^{45}\;, & \omega_3&=e^{34}-e^{52}\;.
\end{align}
Compared to \eqref{eqn:hypo1g}, we have interchanged self\kern.5pt dual and antiself\kern.5pt dual forms; this makes a big
difference since $\omega_3$ is now closed. It can in fact be shown that the \textit{hypo moduli space} (that we have not
considered in the present paper) of $M$ is disconnected, and this explains the existence of this second construction.

In order to solve the new evolution equations, we modify \eqref{eqn:basis} to the new orthonormal basis
\begin{align*}
gh\kern1pt e^1\;,\quad g\kern1pt e^2\;,\quad h\kern1pt e^3\;,\quad g^{-1}e^4\;,\quad h^{-1}e^5\;,\quad dt\;.
\end{align*}
Modifying \eqref{eqn:hypo2g} accordingly, the evolution equations yield
\begin{align}\label{eqn:evol2g}
   2g^2 h\,\partial_tg = 1,\qquad 2gh^2\partial_th = -1.
\end{align}
Observe that $\partial_t(g^2+h^2)=0$. A solution to \eqref{eqn:evol2g} is given by \[g=(1+\sin u)^{\frac12}\;, \qquad
h=(1-\sin u)^{\frac12}\;,\] where \[t=\int_0^u\!\!\cos^2\!x\,dx=\hbox{$\frac12u+\frac14\sin2u$}.\] The resulting metric
is
\[\begin{split}
\cos^2\!u&\>\ee1 + (1+\sin u)\,\ee2 + (1-\sin u)\,\ee3 + \\
&+(1+\sin u)^{-1}\ee4 + (1-\sin u)^{-1}\ee5 + \cos^4\!u\,du\!\otimes du.
 \end{split}\]
  It is defined on
$M\times(-\frac\pi2,\frac\pi2)$, with the change of variable, and has holonomy equal to $\SU(3)$.
Lemma~\ref{eqn:simple3forms} again gives rise to a pencil of Lagrangian submanifolds through each point, and resulting
tri-Lagrangian structure.  \end{example}

\bigskip

\noindent\textbf{Acknowledgements.} This paper overlaps with the first author's \emph{tesi di perfezionamento} at the
Scuola Normale in Pisa, under supervision of the second author; special thanks are due to P. Gauduchon and A.  Swann for
their interest and comments on that work. The authors are also grateful to R.~Bryant, K.~Galicki, N.~Hitchin and
D.~Matessi for various useful conversations that helped the present paper to take shape.

\providecommand{\bysame}{\leavevmode\hbox to3em{\hrulefill}\thinspace}
\providecommand{\MR}{\relax\ifhmode\unskip\space\fi MR }
\providecommand{\MRhref}[2]{%
  \href{http://www.ams.org/mathscinet-getitem?mr=#1}{#2}
}
\providecommand{\href}[2]{#2}

\vskip10pt

\small\noindent Scuola Normale Superiore, Piazza dei Cavalieri 7, 56126 Pisa, Italy\\
\emph{Current address:} Dipartimento di Matematica e Applicazioni, Universit\`a di Milano -- Bicocca, Via Cozzi 53, 20125 Milano, Italy\\
\texttt{diego.conti@unimib.it}\vskip8pt

\noindent Dipartimento di Matematica, Politecnico di Torino, Corso Duca degli Abruzzi 24, 10129 Torino, Italy\\
\texttt{simon.salamon@polito.it}

\begin{thebibliography}{10}

\bibitem{AcharyaDenefHofmanLambert}
B.S. Acharya, F.~Denef, C.~Hofman, and N.~Lambert, \emph{Freund-{Rubin}
  revisited}, hep-th/0308046.

\bibitem{AtiyahWitten}
M.~Atiyah and E.~Witten, \emph{M-theory dynamics on a manifold of {$G_2$}
  holonomy}, Adv. Theor. Math. Phys. \textbf{6} (2003), 1--106.

\bibitem{Bar}
C.~B{\"a}r, \emph{Real {Killing} spinors and holonomy}, Comm. Math. Phys.
  \textbf{154} (1993), no.~3, 509--521.

\bibitem{BarGauduchonMoroianu}
C.~B{\"a}r, P.~Gauduchon, and A.~Moroianu, \emph{Generalized cylinders in
  semi-{R}iemannian and spin geometry}, Math. Z. \textbf{249} (2005), 545--580.

\bibitem{Seminarbericht}
H.~Baum, T.~Friedrich, R.~Grunewald, and I.~Kath, \emph{Twistor and {K}illing
  spinors on {R}iemannian manifolds}, Teubner-Verlag Leipzig/Stuttgart, 1991.

\bibitem{Blair}
D.E. Blair, \emph{{R}iemannian geometry of contact and symplectic manifolds},
  Birkh{\"a}user, 2002.

\bibitem{BoyerGalicki}
C.P. Boyer and K.~Galicki, \emph{3-{S}asakian manifolds}, Surv. Differ. Geom.
  \textbf{7} (1999), 123--184.

\bibitem{Bryant:Calibrated}
R.L. Bryant, \emph{Calibrated embeddings in the special {L}agrangian and
  coassociative cases}, Ann. Global Anal. Geom. \textbf{18} (2000), 405--435.

\bibitem{BryantEtAl}
R.L. Bryant, S.S. Chern, R.B. Gardner, H.L. Goldschmidt, and P.A. Griffiths,
  \emph{Exterior differential systems}, Springer-Verlag, 1991.

\bibitem{Bryant:Private}
R.L. Bryant, \emph{private communication}.

\bibitem{ChiossiFino}
S.~Chiossi and A.~Fino, \emph{Conformally parallel {$G_2$} structures on a
  class of solvmanifolds}, Math. Z. \textbf{252} (2006), no.~4, 825--848.

\bibitem{ChiossiSwann}
S.~Chiossi and A.~Swann, \emph{{$G_2$}-structures with torsion from
  half-integrable nilmanifolds}, J.Geom.Phys. \textbf{54} (2005), 262--285.

\bibitem{thesis}
D.~Conti, \emph{Special holonomy and hypersurfaces}, Ph.D. thesis, Scuola
  Normale Superiore, Pisa, 2005.

\bibitem{Friedrich:OnTheSpinorRepresentationOfSurfaces}
T.~Friedrich, \emph{On the spinor representation of surfaces in {E}uclidean
  3-space}, J.Geom.Phys. \textbf{28} (1998), 143--157.

\bibitem{FriedrichKath}
T.~Friedrich and I.Kath, \emph{Einstein manifolds of dimension five with small
  first eigenvalue of the {D}irac operator}, J. Differential Geom. \textbf{29}
  (1989), 263--279.

\bibitem{FriedrichKim:TheEinsteinDiracEquation}
T.~Friedrich and E.C. Kim, \emph{The {E}instein-{D}irac equation on
  {R}iemannian spin manifolds}, J. Geom. Phys. \textbf{33} (2000), 128--172.

\bibitem{Giovannini}
D.~Giovannini, \emph{Special structures and symplectic geometry}, Ph.D. thesis,
  Universit{\`{a}} degli Studi di {T}orino, 2003.

\bibitem{CalibratedGeometries}
R.~Harvey and H.B. Lawson, \emph{Calibrated geometries}, Acta Math.
  \textbf{148} (1982), 47--157.

\bibitem{Hitchin:StableForms}
N.~Hitchin, \emph{Stable forms and special metrics}, Global Differential
  Geometry: The Mathematical Legacy of Alfred Gray, Contemp. Math., vol. 288,
  American Math. Soc., 2001, pp.~70--89.

\bibitem{Joyce:SpecialLagrangian}
D.~Joyce, \emph{Special {L}agrangian submanifolds with isolated conical
  singularities. {V}. {S}urvey and applications}, J. Differential Geom.
  \textbf{63} (2003), 279--348.

\bibitem{Morel}
B. Morel, \emph{The energy-momentum tensor as a second fundamental form},
  DG/0302205, 2003.

\bibitem{Salamon:Redbook}
S.~Salamon, \emph{Riemannian geometry and holonomy groups}, Pitman Research
  Notes in Mathematics, vol. 201, Longman, Harlow, 1989.

\bibitem{Salamon:ComplexStructures}
\bysame, \emph{Complex structures on nilpotent {L}ie algebras}, J. Pure Appl.
  Algebra \textbf{157} (2001), 311--333.

\end{thebibliography}
\end{document}